\documentclass[11pt]{article}

\usepackage[utf8]{inputenc}
\usepackage{fullpage}
\usepackage{algorithm}
\usepackage{algorithmic}
\usepackage{boxedminipage}
\usepackage{comment}
\usepackage{eqparbox}
\usepackage{multirow}
\usepackage{caption}
\usepackage{subcaption}

\usepackage{amsmath}
\usepackage{amsfonts}
\usepackage{graphicx}

\newlength{\commentindent}
\usepackage{array}
\newcolumntype{R}[1]{>{\raggedleft\let\newline\\\arraybackslash\hspace{0pt}}m{#1}}

\newcommand\BibTeX{{\rmfamily B\kern-.05em \textsc{i\kern-.025em b}\kern-.08em
T\kern-.1667em\lower.7ex\hbox{E}\kern-.125emX}}

\newcommand{\EN}{\mathcal{E}}
\newcommand{\dE}{\partial \EN}

\newcommand{\dive}{\nabla \cdot}
\newcommand{\grad}{\nabla}
\newcommand{\lap}{\Delta}
\newcommand{\jump}[1]{[\![#1]\!]}
\newcommand{\avg}[1]{\{\!\!\{#1\}\!\!\}}

\newcommand{\kernelNl}[2]{\vspace{4pt}\noindent{\bf{\emph{#1 Kernel #2}}}:}

\usepackage[ backend=bibtex, style=authoryear, url=false, isbn=false, maxcitenames=2, maxbibnames=10,giveninits=true,uniquename=init, natbib]{biblatex}
\renewbibmacro{in:}{\ifentrytype{article}{}{\printtext{\bibstring{in}\intitlepunct}}}

\renewbibmacro*{volume+number+eid}{%
  \setunit{\addcomma\space}%
  \printfield{volume}%
  \setunit*{\addnbthinspace}
  \printfield{number}%
  \setunit{\addcomma\space}%
  \printfield{eid}}
\DeclareFieldFormat[article]{number}{\mkbibparens{#1}}

\begin{document}

\title{GPU Acceleration of a High-Order Discontinuous Galerkin Incompressible Flow Solver}

\author{Ali Karakus\thanks{Department of Mathematics, Virginia Tech, 225 Stanger Street, Blacksburg, VA 24061, USA, akarakus@vt.edu},  Noel Chalmers, Kasia \'{S}wirydowicz, \& T.~Warburton}

\maketitle

\begin{abstract}
We present a GPU-accelerated version of a high-order discontinuous Galerkin discretization of the unsteady incompressible Navier–Stokes equations. The equations are discretized in time using a semi-implicit scheme with explicit treatment of the nonlinear term and implicit treatment of the split Stokes operators. The pressure system is solved with a conjugate gradient method together with a fully GPU-accelerated multigrid preconditioner which is designed to minimize memory requirements and to increase overall performance. A semi-Lagrangian subcycling advection algorithm is used to shift the computational load per timestep away from the pressure Poisson solve by allowing larger timestep sizes in exchange for an increased number of advection steps. Numerical results confirm we achieve the design order accuracy in time and space. We optimize the performance of the most time-consuming kernels by tuning the fine-grain parallelism, memory utilization, and maximizing bandwidth. To assess overall performance we present an empirically calibrated roofline performance model for a target GPU to explain the achieved efficiency. We demonstrate that, in the most cases, the kernels used in the solver are close to their empirically predicted roofline performance.
\end{abstract}

\tableofcontents
\section{Introduction}
Finite-element based approximation of the unsteady incompressible Navier-Stokes (INS) equations typically requires high resolution in time and space mandating the use of high performance computing (HPC) techniques. Current trends in HPC show a transition to higher on-node parallelism using accelerators such as Graphical Processing Units (GPUs).  However, developing high-order finite element based flow solvers that take full advantage of modern parallel accelerators is complicated by the need to achieve fine-grain parallelism while effectively exploiting deep non-uniform memory hierarchies. In this work, we focus on the GPU acceleration of a high-order discontinuous Galerkin (DG) spatial disctretization together with semi-implicit temporal discretization combining algebraic splitting and semi-Lagrangian subcycling. 

We choose the discontinuous Galerkin (DG) finite element method for the spatial discretization due to its weak element connectivity and block structured elemental operators. The local stencil of the DG method together with high-order approximations yields highly parallel operators with high arithmetic intensity which are particularly well-suited for GPU accelerators. Kl\"{o}ckner et. al. \citep{klockner_nodal_2009} introduced a GPU
accelerated nodal DG scheme for the first order hyperbolic systems. This approach has since been adapted to, and optimized for, various physical problems \citep{gandham_gpu_2015, modave_gpu_2016,chan_gpu-accelerated_2016,karakus_gpu-accelerated_2016,karakus_gpu_2016}. The implementation and performance optimization of DG methods on GPUs is well documented for first order hyperbolic systems with explicit time integrators. However, only a few papers report similar research regarding optimizing DG discretizations for incompressible flow \citep{roca2011gpu}.

Due to their efficiency for large scale numerical simulations, splitting methods are widely used in time discretizations of the incompressible Navier-Stokes equations. The combination of DG methods with temporal splitting methods has been studied in recent works \citep{ferrer_stability_2014,piatkowski2016stable}. In this work we apply an algebraic splitting technique \citep{chorin1969convergence} as employed in the DG scheme for incompressible flows presented in \citep{shahbazi_high-order_2007}. The reader is referred to \citep{guermond_overview_2006} for an overview of a variety of splitting methods. To further improve the performance of the semi-implicit splitting, we also adopt a semi-Lagrangian subcycling approach, which is closely related to the operator integration factor splitting (OFIS) method \citep{maday_operator_1990}. Stability, dispersion, and dissipation properties of the subcycling approaches are discussed in \citep{giraldo_strong_2003, xiu_strong_2005}. 

Within the algebraic splitting scheme, the velocity and pressure fields are decoupled by enforcing the incompressibility constraint via a Poisson equation for pressure. As we are required to solve this linear system at each time step, preconditioning is applied to overcome the poor conditioning of the Laplacian operator. Multigrid methods \citep{trottenberg_multigrid_2001} are among the most popular and efficient techniques for these equations. Furthermore, a GPU-accelerated version of a unsmoothed aggregation algebraic multigrid (AMG) method \citep{notay2010aggregation} has been investigated recently \citep{gandham_gpu_2014}. However, algebraic multigrid methods require the construction of the full sparse elliptic operator which can lead to high memory requirements. To overcome this limitation, we  use a hybrid multigrid solver as a combination of manually constructed matrix-free $p$-multigrid (pMG) and algebraic multigrid. 


In this work, we present the GPU performance of each of the computationally-intensive kernels present in each step of the temporal splitting scheme. In particular, we show that as more subcycling steps are employed the relative computational cost shifts towards the arithmetically intense non-linear convection kernels. We also show that the majority of the computational costs during the elliptic solvers is contained in the action of the elliptic operators, and we detail the GPU performance of these operators. In order to asses the performance of our computational kernels, we use an empirical roofline model \citep{Volkov2008, Swirydowicz2017}. The model relies on the observation that the GPU is typically a memory-bound device; the runtime of a kernel cannot be faster than the time needed to transfer the data used in the kernel. In addition, the empirical model used in this manuscript takes into account shared memory throughput. Based on the model, we propose a theoretical upper bound for the performance of our code, and this upper bound guides the optimization process. The details of the model are explained in Section 5.

This remainder of this paper is organized as follows. In section 2, we present the mathematical formulation for the DG scheme to approximate the INS equations, including the spatial discretizations and the temporal splitting scheme with semi-Lagrangian approach. Details of the hybrid p-multigrid/ algebraic multigrid solver are given in Section 3, which is followed by numerical validation test cases in Section 4. We then detail key aspects of the GPU implementation, performance analysis and optimization of core kernels in Section 5. Finally, Section 6 is dedicated to concluding remarks and comments on future works. 

\section{Formulation}
We consider a closed two-dimensional domain $\Omega \subset \mathbb{R}^2$ and denote the boundary of $\Omega$ by $\partial\Omega$. We assume that $\partial\Omega$ can be partitioned into two non-overlapping regions, denoted by $\partial\Omega_D$ and $\partial\Omega_N$, along which are prescribed Dirichlet or Neumann boundary conditions, respectively. We are interested in the approximation of the constant density incompressible Navier-Stokes equations
\begin{align}
\frac{\partial \mathbf{u}}{\partial t} + \left(\mathbf{u}\cdot\grad\right)\mathbf{u} &= -\grad p+ \nu\lap \mathbf{u} +\mathbf{f} \label{eq:INS1} \quad &\text{in}\;\;\Omega\times (0,T]\\
\dive\mathbf{u} &= 0 \quad &\text{in}\;\;\Omega\times (0,T],\label{eq:INS2}
 \end{align}
subject to the initial condition 
\begin{equation}\label{eq:initialCondition}
\mathbf{u} = \mathbf{u}_0 \quad \text{for}\;t =0, \mathbf{x} \in \Omega, 
\end{equation}
 and the boundary conditions
\begin{align}\label{eq:DirichletBC}
\mathbf{u} &= \mathbf{g}_D\quad &\text{on}\;\;\mathbf{x}\in\partial\Omega_D, t\in (0,T],\\
\label{eq:NeumannBC} \nu \mathbf{n}\cdot\frac{\partial\mathbf{u}}{\partial \mathbf{x}} - p\mathbf{n} &= 0 \quad &\text{on}\;\;\mathbf{x}\in\partial\Omega_N, t\in (0,T].
\end{align}
Here $\mathbf{u}$ is the velocity field,  $p$ is the static pressure, $\nu$ is the kinematic viscosity, $\mathbf{f}$ is a known body force, and $\mathbf{g}_D$ is prescribed Dirichlet boundary data. In this study, we consider uniform density flows and do not include a density term in the equations above. We discretize this PDE system by first constructing the spatial discretization using the DG method, followed by the temporal discretization  using a temporal splitting scheme. 

\subsection{Spatial Discretization}
We begin by partitioning the computational domain $\Omega$ into $K$ triangular elements $\EN^e$, $e=1,\ldots,K$, such that
\[
\Omega = \bigcup_{e=1}^K \EN^e.
\]
We denote the boundary of the element $\EN^e$ by $\dE^e$. We say that two elements, $\EN^{e+}$ and $\EN^{e-}$, are neighbours if they have a common face, that is $\dE^{e-} \cap \dE^{e+} \neq \emptyset$. We use $\mathbf{n}$ to denote the unit outward normal vector of $\dE$.

We consider a finite element spaces on each element $\EN^e$, denoted $V_N^e  = \mathcal{P}_N(\EN^e)$ where $\mathcal{P}_N(\EN^e)$ is the space of polynomial functions of degree $N$ on element $\EN^e$. As a basis of the finite element spaces we take a set of $N_p = |V_n^e|$ Lagrange polynomials $\{l_n^e\}_{n=0}^{n=N_p}$, interpolating at the Warp $\&$ Blend nodes \citep{warburton_explicit_2006} mapped to the element $\EN^e$. Next, we define the polynomial approximation of the velocity field $\mathbf{u}$ and the pressure field $p$ on each element as 
\begin{align*}
    \mathbf{u}^e &= \sum_{n=0}^{N_p} \mathbf{u}^e_n l_n^e(\mathbf{x}), \\
    p^e &= \sum_{n=0}^{N_p} p^e_n l_n^e(\mathbf{x}),
\end{align*}
for all $\mathbf{x} = (x,y) \in \EN^e$. Using the polynomials $\mathbf{u}^e$ and $p^e$, we introduce the semi-discrete form of the INS system \eqref{eq:INS1}-\eqref{eq:INS2} on an element $\EN^e$ as 
\begin{align}
\label{eq:INS_SD_1}
\frac{d\mathbf{u}^e}{dt} + \mathbf{N}^e(\mathbf{u}^e) &= \mathbf{G}^e p^e + L^e \mathbf{u}^e, \\
D^e \mathbf{u}^e &= 0. \label{eq:INS_SD_1_1}.
\end{align}
Here we have introduced the operators, $\mathbf{N}^e: (V^e_N)^2 \to (V^e_N)^2$, $\mathbf{G}^e: V^e_N \to (V^e_N)^2$, $L^e: V^e_N \to V^e_N$ and $D^e: (V^e_N)^2 \to V^e_N$, which are discrete versions of the nonlinear term $\mathbf{u}\cdot\grad\mathbf{u}$, gradient operator $\grad $, Laplacian $\lap$, and the divergence operator $\dive$, respectively. It remains to define these operators in the DG framework. 

We begin with the discretization of nonlinear term, $\mathbf{u}\cdot\nabla\mathbf{u}$. We use the incompressiblity condition \eqref{eq:INS2} to write $\mathbf{u}\cdot\nabla\mathbf{u}$ in divergence form i.e.,  $\mathbf{u}\cdot\nabla\mathbf{u} = \nabla\cdot\mathbf{F}(\mathbf{u})$, where $\mathbf{F}(\mathbf{u}) = \mathbf{u}\otimes\mathbf{u}$. Multiplying $\mathbf{u}\cdot\nabla\mathbf{u}$ by a test function $v \in V^e_N $, integrating over the element $\EN^e$, and performing integration by parts, we define the discrete nonlinear term $\mathbf{N}^e(\mathbf{u})$ via the following variational form
\begin{equation}
\label{eq:INS_SD_2}
(v,\mathbf{N}^e(\mathbf{u}^e))_{\EN^e} =
-(\nabla v,\mathbf{F}(\mathbf{u}^e))_{\EN^e} 
+ (v,\mathbf{n}\cdot \mathbf{F}^{*})_{\dE^{e}}
\end{equation}
Here we have introduced the inner product $(u , v)_{\EN^e} $ to denote the integration of the product of $u$ and $v$ computed over the element $\EN^e$ and, analogously,  the inner product $(u , v)_{\dE^e}$ to denote the integration along the element boundary $\dE^e$. 

Due to the discontinuous approximation space, the flux function $\mathbf{F}$ is not uniquely defined in the boundary inner product and hence, it is replaced by a numerical flux function $\mathbf{F}^*$ which depends on the local and neighboring traces values of $\mathbf{u}$ along $\partial\EN^e$. One each element we denote the local trace values of $\mathbf{u}^e$ as $\mathbf{u}^-$ and the corresponding neighboring trace values as $\mathbf{u}^+$. Note that we will suppress the use of the $e$ superscript when it is clear which element is the local trace. Using this notation we choose as a numerical flux $\mathbf{F}^*$ the local Lax-Friedrichs numerical flux, i.e.,
\begin{eqnarray}
\label{eq:INS_SD_3}
\mathbf{F}^{*}=   \avg{\mathbf{F}(\mathbf{u})}  + \frac{1}{2} \mathbf{n} \Lambda^e \jump{\mathbf{u}}.
\end{eqnarray}
Here we use the notation $\avg{\mathbf{u}}$ and $\jump{\mathbf{u}}$ to denote the average and jump of $\mathbf{u}$ along the the trace $\dE^e$, that is
\begin{equation}\label{Eq.AverageJumpScalar}
\avg{\mathbf{u}} = \frac{\mathbf{u}^{+} + \mathbf{u}^{-}}{2}, \quad  \jump{\mathbf{u}^e} = \mathbf{u}^{+} -\mathbf{u}^{-}.
\end{equation}
The parameter $ \Lambda$ in \eqref{eq:INS_SD_3} is a stabilization parameter, which introduces artificial diffusion required to stabilize the numerical discretization of the nonlinear term. The parameter is chosen to be the maximum eigenvalue of the flux Jacobian in absolute value, i.e.  
\[
\Lambda = \max_{\mathbf{u}\in[\mathbf{u}^-,\mathbf{u}^+]} \left| \mathbf{n}\cdot\frac{\partial \mathbf{F} }{\partial \mathbf{u}}\right|.
\]
The choice of local Lax-Friedrichs flux leads to a stable and easily evaluated numerical flux function. In the case of Dirichlet boundaries $\dE^e \cap \Omega_D \neq \emptyset$, we weakly enforce the Dirichlet boundary condition \eqref{eq:DirichletBC} by choosing $\mathbf{u}^{+} = \mathbf{g}_D$ along this trace, while for Neumann boundaries $\dE^e \cap \Omega_N \neq \emptyset$, we simply choose $\mathbf{u}^{+} = \mathbf{u}^{-}$.

Moving on to the gradient and divergence operators, $\mathbf{G}^e$ and $D^e$, respectively, we use the DG approximation to discretize these operators in a way analogous to that described above for the nonlinear operator $\mathbf{N}^e(\mathbf{u})$. Namely, we multiply the pressure gradient $\nabla p^e$ and the velocity divergence $\nabla\cdot\mathbf{u}^e$ by a test function $v\in V^e_N$, integrate over the element $\EN^e$, and integrate by parts twice. We choose the numerical fluxes $p^*$ and $\mathbf{u}^*$ to be simply the central fluxes $p^* = \avg{p}$ and $\mathbf{u}^*= \avg{\mathbf{u}}$ to obtain the following variational definitions of $\mathbf{G}^e$ and $D^e$
\begin{align}
\label{eq:INS_SD_4_1}(v,\mathbf{G}^e p^e)_{\EN^e} &= (v,\grad p^e)_{\EN^e} 
    + \frac{1}{2}(v,\mathbf{n}\jump{p})_{\dE^{e}}, \\
\label{eq:INS_SD_4_2}        (v,D^e \mathbf{u}^e)_{\EN^e} &= (v,\dive\mathbf{u}^e)_{\EN^e} + \frac{1}{2}(v,\mathbf{n}\cdot \jump{\mathbf{u}})_{\dE^{e}}.
\end{align}
We impose boundary conditions for these operators slightly differently than for $\mathbf{N}^e(\mathbf{u})$. Specifically, along Dirichlet boundaries we take $\mathbf{u}^{*} = \mathbf{g}_D$ and $p^* = p^{-}$ and for Neumann boundaries, we choose $\mathbf{u}^{*} = \mathbf{u}^{-}$ and $p^* = 0$.

Finally, to discretize the Laplacian operator $L^e$, we note that $\lap\mathbf{u} = \dive\grad \mathbf{u}$ holds in the continuous setting. The Laplacian operators can then be discretized for a DG method by simply using the composition of the discrete gradient and divergence operators so that $L^e = D^e\cdot\mathbf{G}^e$. This leads to the well-known local DG discretization. Forming the above-mentioned composition and applying integration by parts to the volume term leads to the following variational definition of the Laplacian operator $L^e$ 
\begin{equation*}
    (v,L^e\mathbf{u}^e)_{\EN^e} = -(\grad v ,\grad\mathbf{u}^e)_{\EN^e} + (v,\mathbf{n}\cdot\grad\mathbf{u}^*)_{\dE^{e}} - (\mathbf{n}\cdot\grad v,\mathbf{u}^*-\mathbf{u}^{-})_{\dE^{e}}.
\end{equation*}
In contrast to the gradient and divergence operators, simply choosing central fluxes for $\mathbf{n}\cdot\grad\mathbf{u}^*$ and $\mathbf{u}^*$ results in an inconsistent and weakly unstable scheme \citep{zhang_analysis_2003}. We therefore follow the Symmetric Interior Penalty DG (SIPDG) approach  \citep{wheeler1978elliptic,arnold_interior_1982} and choose the numerical flux terms to be the central fluxes augmented by the penalty term, i.e., $\mathbf{u}^* = \avg{\mathbf{u}}$ and $\mathbf{n}\cdot\grad\mathbf{u}^* = \mathbf{n}\cdot \avg{\grad\mathbf{u}} + \tau \jump{\mathbf{u}}$. The variational form can then be written 
\begin{align}\label{eq:INS_SD_5}
    (v,L^e\mathbf{u}^e)_{\EN^e} = -&(\grad v ,\grad\mathbf{u}^e)_{\EN^e} + (v,\mathbf{n}\cdot\avg{\grad\mathbf{u}})_{\dE^{e}}  \\ \nonumber &-\frac{1}{2}(\mathbf{n}\cdot\grad v,\jump{\mathbf{u}})_{\dE^{e}} + (v,\tau\jump{\mathbf{u}})_{\dE^{e}}.
\end{align}
The penalty parameter $\tau$ must be chosen to be sufficiently large in order to enforce coercivity. Care must be taken, however, as selecting large $\tau$ results in poor conditioning of the Laplacian operator and degrades the performance of linear solvers. Along each face $\dE^{ef} = \EN^{e+} \cap \EN^{e-}$, we select a penalty parameter $\tau^{ef}$ using the lower bound estimate derived in \citep{shahbazi_explicit_2005}: 
\begin{equation}
\label{Eq.Ch2.PenaltyParameter}
\tau^{ef}=\frac{(N+1)(N+2)}{2}\max\left(\frac{1}{h^{ef}_+}, \frac{1}{h^{ef}_-}\right),
\end{equation}
where $h^{ef}_+$ and $ h^{ef}_-$ are characteristic length scales of the elements $\EN^{e+}$ and $\EN^{e-}$ on either side of the face $\dE^{ef}$ and are defined as $h^{ef}_+ = \frac{|\EN^{e+}|}{|\dE^{ef}|}$ and $h^{ef}_- = \frac{|\EN^{e-}|}{|\dE^{ef}|}$. Once the penalty parameter is chosen large enough to enforce coercivity the SIPDG discretization gives a high-order accurate discretization of the Laplacian operator. Boundary conditions for the discretized Laplacian operator are imposed in a way analogous to that described for the gradient and divergence operators above.

System \eqref{eq:INS1} together with the definitions of discrete operators $\mathbf{N}^e$, $L^e$, $\mathbf{G}^e$ and $D^e$ in \eqref{eq:INS_SD_2}, \eqref{eq:INS_SD_4_1}, \eqref{eq:INS_SD_4_2} and \eqref{eq:INS_SD_5} completes the semi-discrete form of the scheme in \eqref{eq:INS_SD_1}. In the next section, we proceed to the fully discrete scheme by introducing the semi-explicit time integration method and semi-Lagrangian subcycling approach.  

\subsection{Temporal Discretization}
Assembling the semi-discrete system in \eqref{eq:INS_SD_1} on each element $\EN$ into global system, we arrive to following global problem
\begin{subequations}
\label{eq:INS_TD_1}
\begin{equation}
\frac{\partial \mathbf{U}}{\partial t} + \mathbf{N(U)} = -\mathbf{G} P+\mathbf{L U},
\end{equation}
\begin{equation}
D\mathbf{U} = 0.
\end{equation}
\end{subequations}
To simplify the notation, we use capital letters and drop the superscript $e$ to denote the global assembled vectors of the degrees of freedom. 

We implement a high-order temporal discretization of the flow equations by adopting an $S$ order backward differentiation method for the stiff diffusive term $\mathbf{L U}$ and an $S$ order extrapolation method for non-linear advective term $\mathbf{N(U)}$. With this formulation, \eqref{eq:INS_TD_1} can be advanced from time level $ t^n $ to $ t^{n+1} = t^{n} + \Delta t$ by solving the equation, 
\begin{subequations}
\label{eq:INS_TD_2}
    \begin{equation}
        \label{eq:discrete-poly}
        \gamma \mathbf{U}^{n+1} = \sum_{i=0}^S \beta_i \mathbf{U}^{n-i} -\Delta t \sum_{i=0}^S \alpha_i \mathbf{N}(\mathbf{U}^{n-i}) + \nu \Delta t L \mathbf{U}^{n+1} - \Delta t \mathbf{G} P^{n+1},
    \end{equation}
    \begin{equation}
        \label{eq:discrete-incompressibility}
        \mathbf{D}\cdot\mathbf{U}^{n+1} = 0.
    \end{equation}
\end{subequations}
where the coefficients $ \beta $, and $ \gamma$ correspond to the stiffly stable backwards differentiation scheme and the coefficients $ \alpha$ correspond to the extrapolation scheme. For the second order scheme the coefficients are $\gamma = 3/2 $, $\beta_0= 2 $, $\beta_1 = 1/2 $ and $\alpha_0= 2$, $\alpha_1 = -1 $. Because this high-order explicit evaluation is not self starting, it is initialized with lower order counterparts; their values can be found in \citep{karniadakis_spectral/hp_2005}.  

We replace the fully discrete scheme \eqref{eq:INS_TD_2} with an algebraically split version following \citep{shahbazi_high-order_2007} in order to solve for velocity and pressure separately instead of solving a fully coupled system. To do this, we first introduce $\delta^k P^{n+1}$ to denote the high-order backward finite differences of pressure, defined recursively as
$\delta^k P^{n+1} = \delta^{k-1} P^{n+1} - \delta^{k-1} P^{n}$
and $\delta^0 P^{n} = P^{n}$. We also introduce the difference $\sigma^k P^n = P^{n+1} - \delta^k P^{n+1}$ where $\sigma^k P^n$ does not depend on $P^{n+1}$. Using this notation, algebraic splitting scheme can be written in four steps as follows, 
\begin{subequations}
\label{eq:INS_TD_3}
    \begin{equation} 
        \label{eq:INS_TD_3_1} 
        \mathbf{\hat{U}} = \sum_{i=0}^S \beta_i \mathbf{U}^{n-i} - \Delta t \sum_{i=0}^S \alpha_i \mathbf{N}(\mathbf{U}^{n-i}).
    \end{equation}
    \begin{equation} 
        \label{eq:INS_TD_3_2} 
        \left(-L + \frac{\gamma}{\nu\Delta t}\mathcal{I}\right)\hat{\hat{\mathbf{U}}} = \frac{1}{\nu\Delta t}\hat{\mathbf{U}} -  \frac{1}{\nu} \mathbf{G} \sigma^{S+1} P^{n}.
    \end{equation}
    \begin{equation} 
        \label{eq:INS_TD_3_3} 
        -L \delta^{S+1} P^{n+1} = -\frac{\gamma}{\Delta t}\mathbf{D}\cdot\hat{\hat{\mathbf{U}}}.
    \end{equation}
    \begin{equation}
    \label{eq:INS_TD_3_4}
    \begin{aligned}
        \mathbf{U}^{n+1} &= \hat{\hat{\mathbf{U}}} - \frac{\Delta t}{\gamma} \mathbf{G} \delta^{S+1} P^{n+1}, \\
        P^{n+1} &\leftarrow \delta^{S+1} P^{n+1} + \sigma^{S+1} P^{n}.
        \end{aligned}
    \end{equation}
\end{subequations}
The steps of this splitting scheme can be interpreted as 1) a pure advection evaluation in \eqref{eq:INS_TD_3_1}, 2) a screened Poisson equation in \eqref{eq:INS_TD_3_2} to implicitly step the diffusive term, 3) a pressure correction in \eqref{eq:INS_TD_3_3} to enforce divergence free velocity, and finally 4) a corrective update step in \eqref{eq:INS_TD_3_4}. This splitting scheme reduces the cost of the temporal discretization to a combination of explicit steps and two linear elliptic solves. The maximum stable time step size will still be determined by the spectrum of the convective term $\mathbf{N}(\mathbf{U})$ and the elliptic solves will still dominate the cost of each time step. To reduce the computational cost of each time step we consider a subcycling method to increase the size of the maximum stable time step. 

\subsection{A Lagrangian Subcycling Method}

The stable timestep size of the splitting scheme is restricted by a Courant-Friedrichs-Lewy (CFL) condition as a result of the explicit treatment of the convective term $\mathbf{N}(\mathbf{U})$. To overcome this restriction, we implement a semi-Lagrangian subcycling method for the INS equations which can be viewed as a high-order operator integration factor splitting approach of Maday et.al. \citep{maday_operator_1990}, and is similar to the semi-Lagrangian subcycling approach presented in \citep{xiu_strong_2005}.

The splitting scheme \eqref{eq:INS_TD_3_1}-\eqref{eq:INS_TD_3_4} provides a natural setting for the subcycling method by separating the advection step from the elliptic parts. We consider the explicit advective stage \eqref{eq:INS_TD_3_1} which approximates an explicit time step of the total derivative $\frac{D\mathbf{U}}{Dt} \equiv \frac{\partial \mathbf{U}}{\partial t} + \mathbf{U}\cdot\nabla\mathbf{U}$. In the Lagrangian frame, we can replace this stage with
\begin{equation} 
        \mathbf{\hat{U}} = \sum_{i=0}^S \beta_i \mathbf{\tilde{U}}^{n-i}.
\end{equation}
where $\tilde{\mathbf{U}}^n$ is the Lagrangian velocity field at time $t^n$. Since, in our time stepping scheme we hold only the history of the velocity fields in the Eulerian frame, i.e. $\mathbf{U}^{n-i}$ for $i=0,\ldots,S$, it remains to show how to compute the Lagrangian velocities from the Eulerian history. 

As described in \citep{maday_operator_1990} and \citep{xiu_strong_2005} the Lagrangian velocity field $\tilde{\mathbf{U}}^{n-i}$ can be approximated by time-stepping the following subproblem
\begin{equation}
  \begin{aligned}
	\label{eq:SS_2}
     \frac{\partial \tilde{\mathbf{U}}_i} {\partial t} &= - \bar{\mathbf{U}}\cdot \nabla\tilde{\mathbf{U}}_i,  \\
     \tilde{\mathbf{U}}_i\left(\mathbf{x},t^{n-i}\right) &= \mathbf{U}^{n-i}\left(\mathbf{x}\right),
\end{aligned}    
\end{equation}
from $t^{n-i}$ to $t^{n+1}$ and setting $\tilde{\mathbf{U}}^{n-i} = \tilde{\mathbf{U}}^i(\mathbf{x},t^{n+1})$. Here the advective velocity field $\bar{\mathbf{U}}(\mathbf{x},t)$ is a degree $S$ polynomial in $t$ interpolating the Eulerian velocities $\mathbf{U}^{n-i}$ at $t=t^{n-i}$ for $i=0,\ldots,S$, respectively. 

Discretizing the linear system \eqref{eq:SS_2} using the DG formulation on each element $\EN^e$ by an analogous procedure to that used above we obtain the semi-discrete system
\begin{equation}
	\label{eq:SS_3}
     \frac{\partial \tilde{\mathbf{U}}^e_i} {\partial t} = - \tilde{\mathbf{N}}^e(\bar{\mathbf{U}^e},\tilde{\mathbf{U}}^e_i),
\end{equation}
where the operator $\tilde{\mathbf{N}}^e(\bar{\mathbf{U}}^e,\tilde{\mathbf{U}}^e_i)$ is defined as satisfying the following variational statement
\begin{equation}\label{eq:INS_CUB_N}
    (v,\tilde{\mathbf{N}}^e(\bar{\mathbf{U}}^e,\tilde{\mathbf{U}}^e_i))_{\EN^e} = -(\nabla v,\tilde{\mathbf{F}}(\bar{\mathbf{U}}^e,\tilde{\mathbf{U}}^e_i))_{\EN^e} + (v,\mathbf{n}\cdot \tilde{\mathbf{F}}^{*})_{\dE^{e}},    
\end{equation}
for all $v \in V^e_N$. Here $\tilde{\mathbf{F}}(\bar{\mathbf{U}}^e,\tilde{\mathbf{U}}^e_i) = \bar{\mathbf{U}}^e \otimes \tilde{\mathbf{U}}^e_i$ and we have used the fact the Eulerian velocity fields are divergence-free in order to write $\bar{\mathbf{U}}\cdot \nabla\tilde{\mathbf{U}}_i = \nabla\cdot\tilde{\mathbf{F}}(\bar{\mathbf{U}}^e,\tilde{\mathbf{U}}^e_i)$. We again choose the local Lax-Friedrichs flux in the definition of $\tilde{\mathbf{N}}^e(\bar{\mathbf{U}^e},\tilde{\mathbf{U}}^e_i)$, i.e. we take
\[
\tilde{\mathbf{F}}^{*}= \avg{\tilde{\mathbf{F}}(\bar{\mathbf{U}}^e,\tilde{\mathbf{U}}^e_i)}  + \frac{1}{2} \mathbf{n} \tilde{\Lambda} \jump{\tilde{\mathbf{U}}^e_i},
\]
where 
\[
\tilde{\Lambda} = \max_{\tilde{\mathbf{U}}\in[\tilde{\mathbf{U}}^-_i,\tilde{\mathbf{U}}^+_i]} \left| \mathbf{n}\cdot\frac{\partial \tilde{\mathbf{F}} }{\partial \tilde{\mathbf{U}}} \right|.
\]
We can compute this operator by splitting its evaluation into volume and surface integral contributions. 

We time step each of the subproblems \eqref{eq:SS_3} for $i=0,\ldots,S$ with a fourth-order low-storage explicit Runge-Kutta (LSERK) method  \citep{williamson1980low,carpenter1994fourth}. We denote by $\Delta t_s $ the timestep size used in this LSERK scheme and take the the macro timestep size $\Delta t$ to be a multiple of $\Delta t_s$, i.e.  $\Delta t = N_s \Delta t_s$. In this way, we say that we use $N_s$ advection subcycles per time step of the full INS system. 

Since the CFL condition now only limits the size of the LSERK timestep $\Delta t_s$ this subcycling approach enables using $N_s$ times larger macro timesteps, hence $N_s$ times fewer linear solves, per macro time step. We instead require $S\times N_s$ additional explicit advection steps using the linearity of \eqref{eq:SS_2} in $\bar{\mathbf{U}}$ and applying superposition. The efficiency of the subscycling method therefore comes from the fast evaluation of these advection steps using the DG discretization which does not require global mass matrix inversion. Note, however, that increasing the macro timestep size effects the performance of screened Poisson solve in \eqref{eq:INS_TD_3_2}. In Section 4, we briefly discuss the benefit of the subcycling method on the total solver time, and the impact on the performance of the screen Poisson equation solver.

\section{Linear Solvers}

Each time step of the temporal splitting discretization \eqref{eq:INS_TD_3} requires solving discrete screened Poisson equation \eqref{eq:INS_TD_3_2} and discrete Poisson problem \eqref{eq:INS_TD_3_3}. We must therefore ensure that these linear systems are solved as fast and as efficiently as possible. For large meshes and/or high degree $N$, assembling a full matrix and using a direct solver is not feasible. Thus, we resort to iterative solvers and, noting that the IP discretization \eqref{eq:INS_SD_5} is symmetric positive-definite with our chosen penalty parameter, we choose a preconditioned conjugate gradient (PCG) iterative method to solve \eqref{eq:INS_TD_3_2} and \eqref{eq:INS_TD_3_3}. 

For the screened Poisson problem in \eqref{eq:INS_TD_3_2}, we note that since the time step $\Delta t$ is usually small, the screened Poisson operator is dominated by the mass matrix with coefficient $\frac{1}{\nu \Delta t}$. Since the mass matrix is block diagonal and the elemental geometric factors are constant on each triangular/tetrahedral element, this mass matrix operator is simple and inexpensive to invert. We therefore choose the scaled inverse mass matrix on each element as a preconditioner for the screened Poisson problem \eqref{eq:INS_TD_3_2}. As we detail below, this preconditioner is usually an effective choice, however, the number of PCG iterations required to solve \eqref{eq:INS_TD_3_2} increases when the number of subcycling steps is increased due to a larger time step size $\Delta t$. 

For the Poisson problem in \eqref{eq:INS_TD_3_3}, we consider two types of multigrid preconditioners. The first is a purely algebraic multigrid (AMG) preconditioner \citep{stuben2001review}. The coarse levels of this AMG method are constructed as unsmoothed aggregations of maximal independent node sets, see \citep{notay2006aggregation, notay2010aggregation}, while smoothing is chosen to be a degree 2 Chebyshev iteration \citep{adams2003parallel}. The multigrid preconditioning cycle itself consists of a K-cycle on the finest two levels, followed by a V-cycle for the remaining coarse levels. We choose these components of the AMG preconditioner to obtain, as presented in \citep{gandham_gpu_2014}, a fine-grain parallel multigrid operation, i.e., the sparse stiffness matrix, the sparse prolongation and restriction actions, and the smoothing operations are all simple to parallelize on the GPU. 

The PCG method using this full AMG preconditioner performs reasonably well but the iteration counts do scale roughly linearly with degree $N$. Furthermore, a significant amount of storage is required to construct a full stiffness matrix for higher degrees. Hence, we consider a multigrid preconditioner where we manually coarsen from degree $N$ to degree 1 before setting up the same AMG levels for the degree 1 coarse stiffness matrix. This approach is similar to that considered in \citep{lottes2005hybrid}, which combined Schwartz patch smoothers on manually constructed degree $p$ multigrid levels before proceeding to a degree 1 coarse problem. With this manual coarsening approach, we are able to implement the finest levels of the multigrid cycle in a matrix-free way and avoid the storage of the full degree $N$ stiffness matrix. We refer to this hybrid manual/algebraic multigrid preconditioner as pMG-AMG.  

\section{Numerical Tests}
In this section we present two dimensional benchmark tests to verify the spatial and temporal accuracy of the proposed scheme and show the performance of the pMG-AMG and AMG preconditioners for the Poisson solver. We then continue with the flow past a square cylinder test problem to describe relative importance of each solver step in the splitting scheme. We also show the effects of using semi-Lagrangian subcycling on relative runtimes and on performance of implicit solves. This test case will inform our later discussion regarding GPU implementations and kernel optimization discussed in the next section.

In all the test cases, unless explicitly stated otherwise, we use the second-order time splitting scheme i.e. we use second-order backward differentiation and extrapolation and use the first-order pressure increment.

\subsection{Taylor Vortex}
Taylor vortex problem is used to test the temporal and spatial accuracy of the method. The solution is known everywhere for all times and given by
\begin{equation}
 \begin{aligned}
\label{Eq.Taylor_Eaxact}
\mathbf{u} &= \left(-\sin(2\pi y)e^{-\nu 4 \pi^2 t}\right)\mathbf{i} + \left(\sin(2\pi x)e^{-\nu 4 \pi^2 t}\right)\mathbf{j}\\
p &= -\cos(2\pi x)\cos(2 \pi y)e^{-\nu 8 \pi^2 t}.
 \end{aligned}
\end{equation}
This flow test is performed with $\nu=0.01$ and is run until the final time $T=3$ is reached at which point the velocity field decays to approximately one-third of its initial amplitude. The computational domain of $[-0.5,0.5]^2$ is discretized with a mesh of unstructured triangular elements. The domain boundaries are specified to be inflow boundaries at the upper, lower and left walls while the right wall is specified to be an outflow boundary. We specify the exact Dirichlet boundary condition for velocity/pressure at the inflow/outflow boundaries, respectively. 

\begin{figure}[htp]
\begin{center}
 \begin{subfigure}{0.45\textwidth}
   \includegraphics[width=0.99\textwidth]{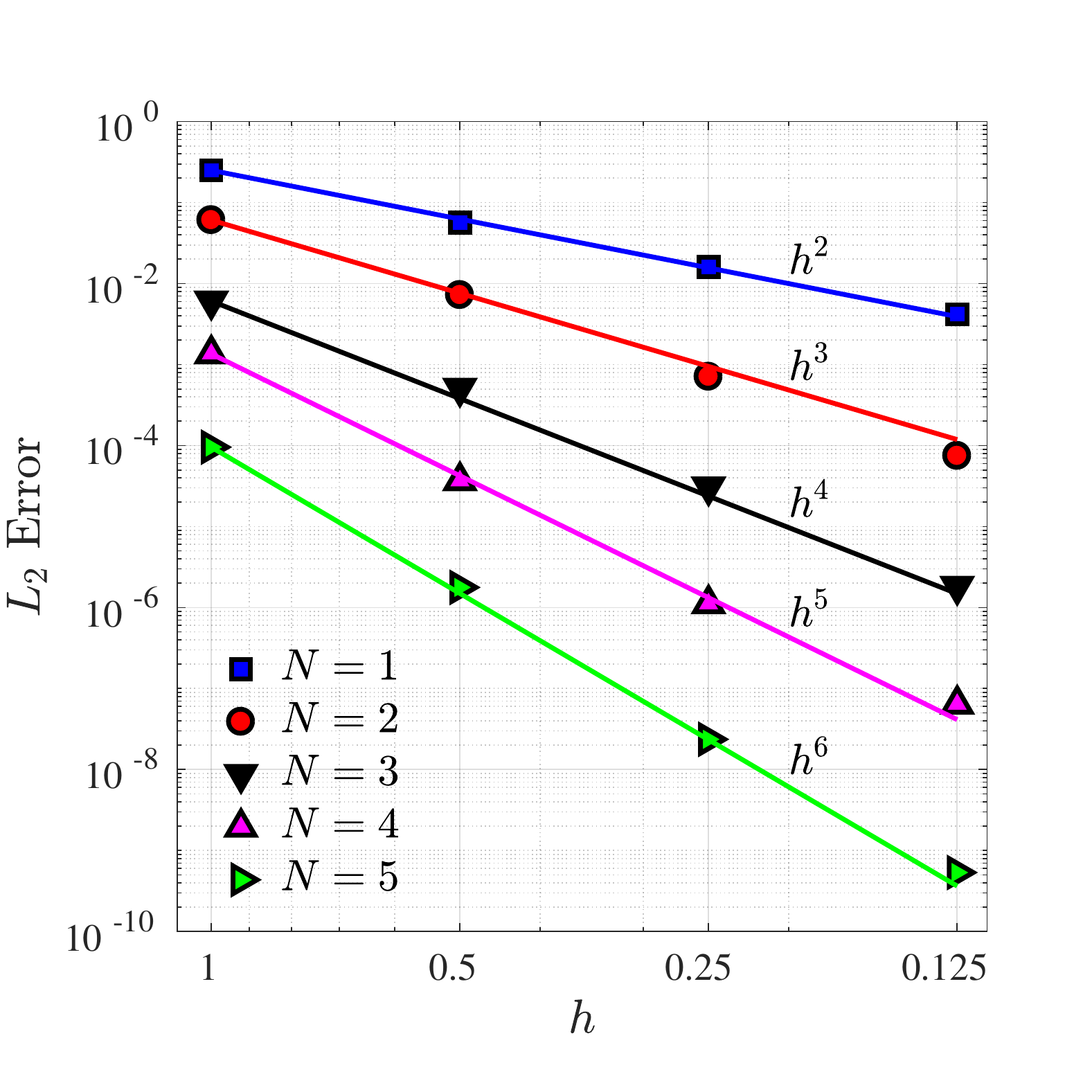}
 \end{subfigure}
 \begin{subfigure}{0.45\textwidth}
   \includegraphics[width=0.99\textwidth]{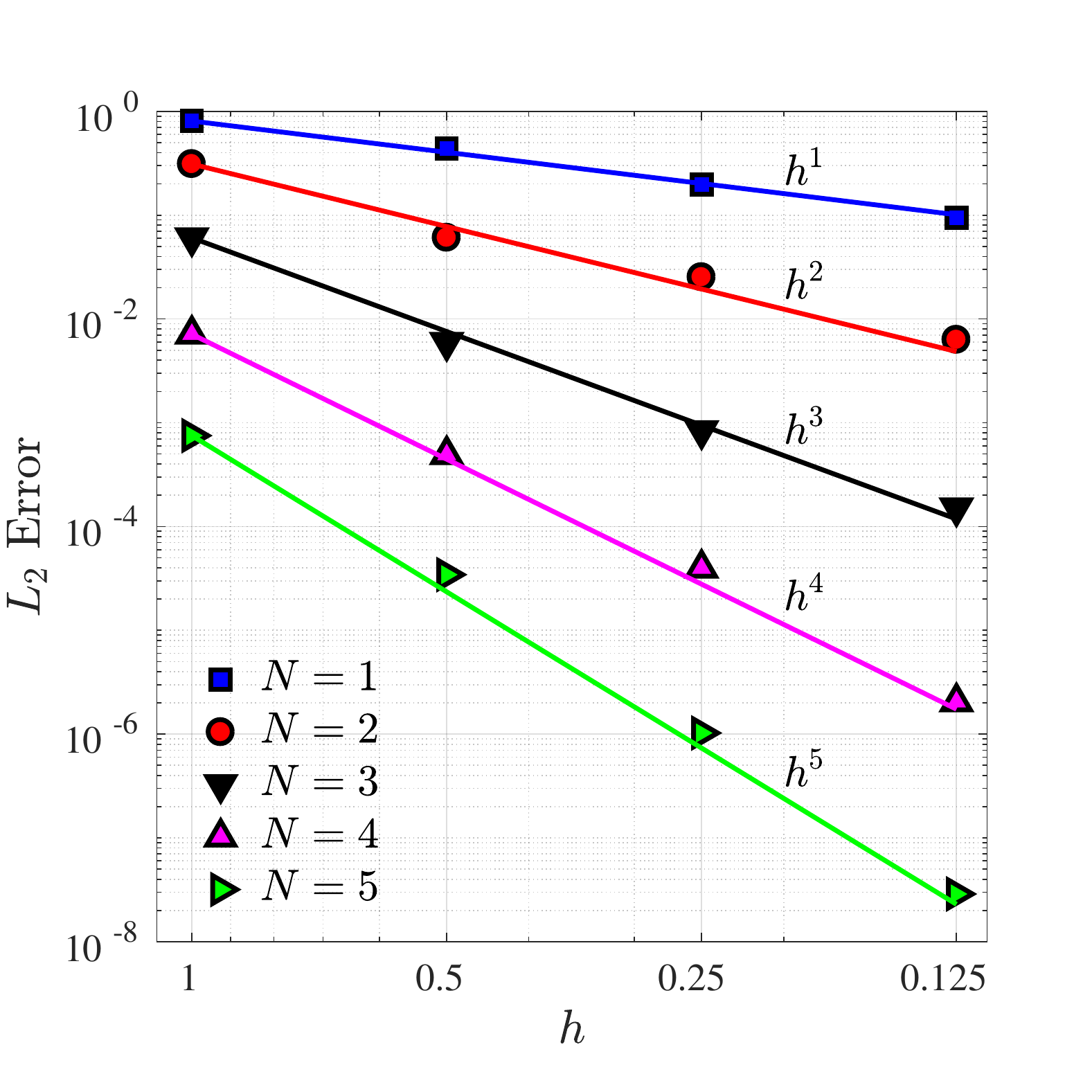}
 \end{subfigure}
\caption{Spatial accuracy test for the Taylor vortex test problem \eqref{Eq.Taylor_Eaxact} using $L_2$ relative errors on successively refined triangular elements. The error in the $x$-velocity is shown on the left and the error in the pressure is shown right.}
\label{fig:SpatialConvergence2D}
\end{center}
\end{figure}

Figure \ref{fig:SpatialConvergence2D} shows the computed $L_2$ norm of the numerical error in the pressure and the $x$ component of velocity at the final time $T=3$. We begin with an unstructured mesh of $K=35$ elements and carry out convergence study with successive $h$ refinement and several degrees $N$. The figure demonstrates the expected $h^{N+1}$ and $h^{N}$ convergence rate in the numerical error. The $y$-velocity has similar convergence properties as the $x$-velocity and is not shown in the figure.  
\begin{figure}[htb!]
\begin{center}
 \begin{subfigure}{0.45\textwidth}
   \includegraphics[width=0.99\textwidth]{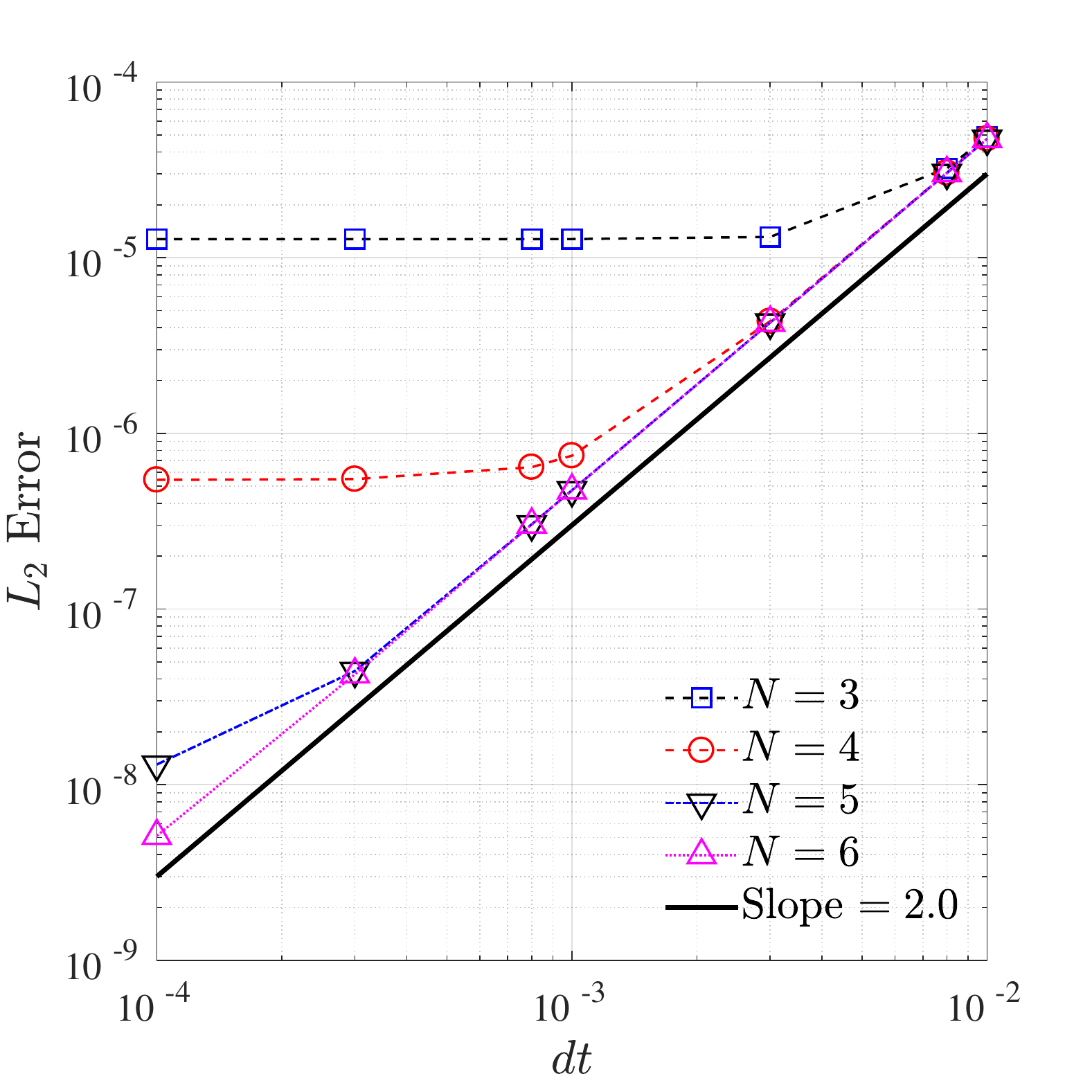}
   \caption{}
   \label{fig:VortexTemporal_1}
 \end{subfigure}
 \begin{subfigure}{0.45\textwidth}
   \includegraphics[width=0.99\textwidth]{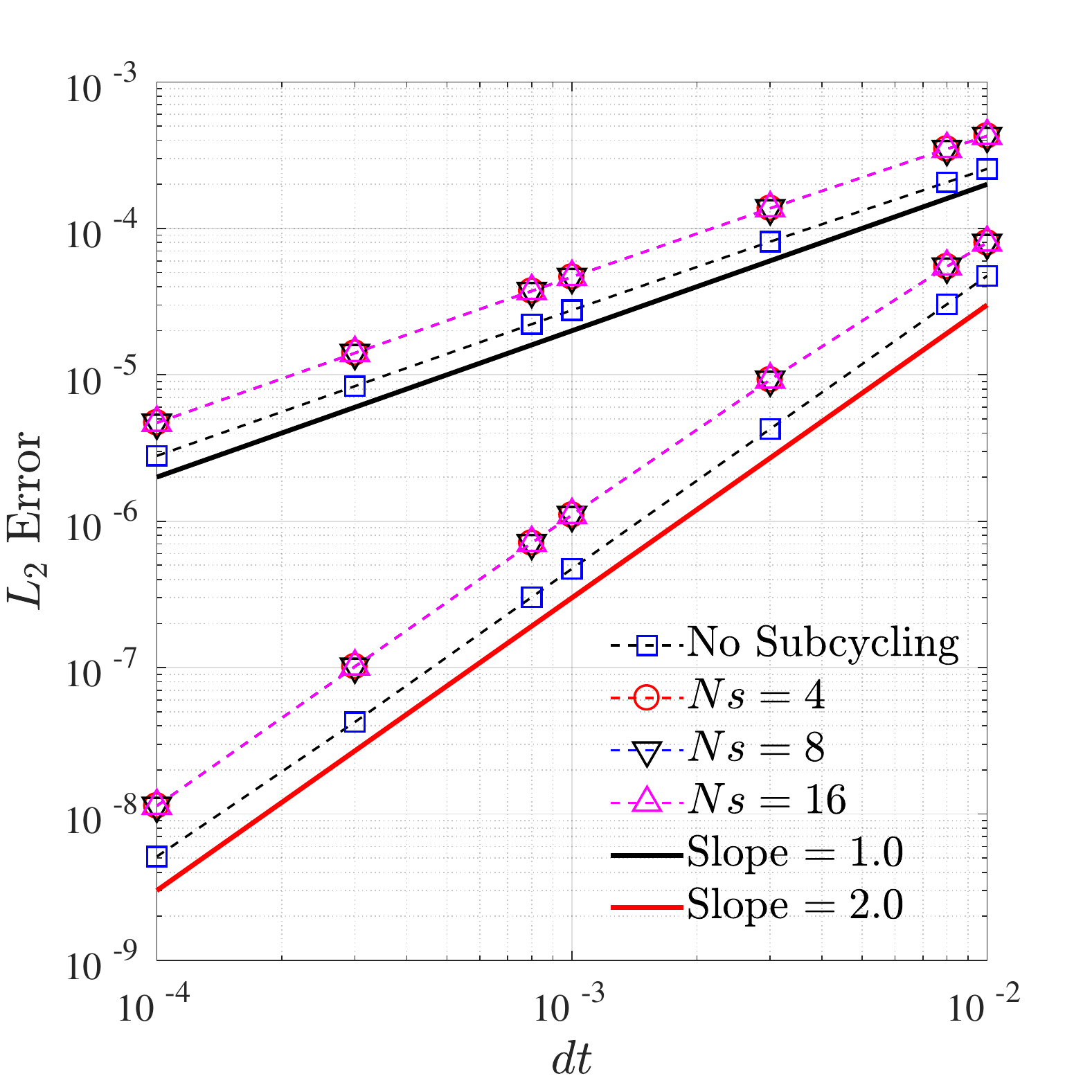}
   \caption{}
   \label{fig:VortexTemporal_2}
 \end{subfigure}
 \begin{subfigure}{0.45\textwidth}
   \includegraphics[width=0.99\textwidth]{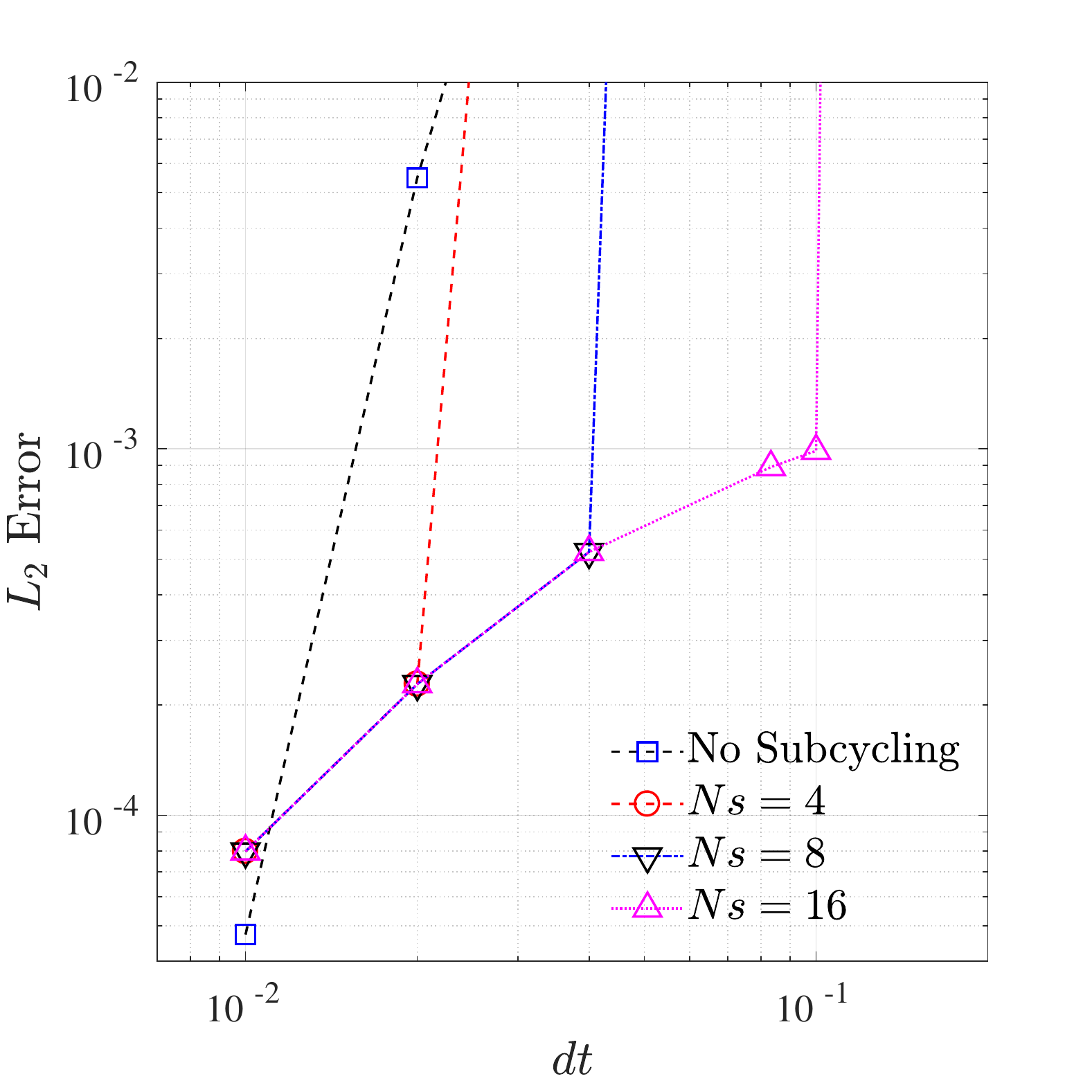}
   \caption{}
   \label{fig:VortexTemporal_3}
 \end{subfigure}
\caption{Temporal accuracy test for Taylor vortex test problem using $L_2$ relative errors of $x$-velocity (a)  timestep refinement study for different orders of approximation. (b) Comparison of first and second order subcycling approaches in stable timestep region for $N=6$. (c) Maximum stable timestep size for  different number of substeps for $N=6$.}
\label{fig:VortexTemporal}
\end{center}
\end{figure}

In Figure \ref{fig:VortexTemporal_1} we show the $L_2$ error of the $x$-velocity in a timestep refinement study. For low-order approximations spatial error dominates the temporal error and decreasing the time step size further does not improve the accuracy. The expected second order accuracy is obtained for all the cases in the region where the temporal errors dominate. The pressure and the $y$-velocity exhibit similar temporal convergence properties and are not included. 

We show in Figure \ref{fig:VortexTemporal_2} the $L_2$ norm of the relative error for the $x$-velocity for subcycling with different number of substeps and without subcycling with the first and the second order time integration and $N=6$. Although, there is no computational advantage of using subscycling if the time step size is stable for standard integration, we include the figure to show the formal accuracy of the method. Subcycling shows the expected first and second order accuracy that is independent from the number of substeps. The numerical error depends on the macro timestep size, $dt$ for the problems with the same spatial resolutions. Comparing with temporal integration without subcycling, we observe slightly larger errors in the subcycling approach. This shift in the error can be explained by the dissipation added to the scheme to stabilize the system with high CFL numbers. Finally, the $L_2$ norm of the numerical error for the subcycling method with varying number of substeps is shown in Figure \ref{fig:VortexTemporal_3} for $N=6$ as the timestep size increases. We see in this figure that the numerical error remains controlled for larger time steps sizes as we take more subcycling steps.     

\begin{figure}[htb!]
\begin{center}
 \begin{subfigure}{0.45\textwidth}
   \includegraphics[width=0.99\textwidth]{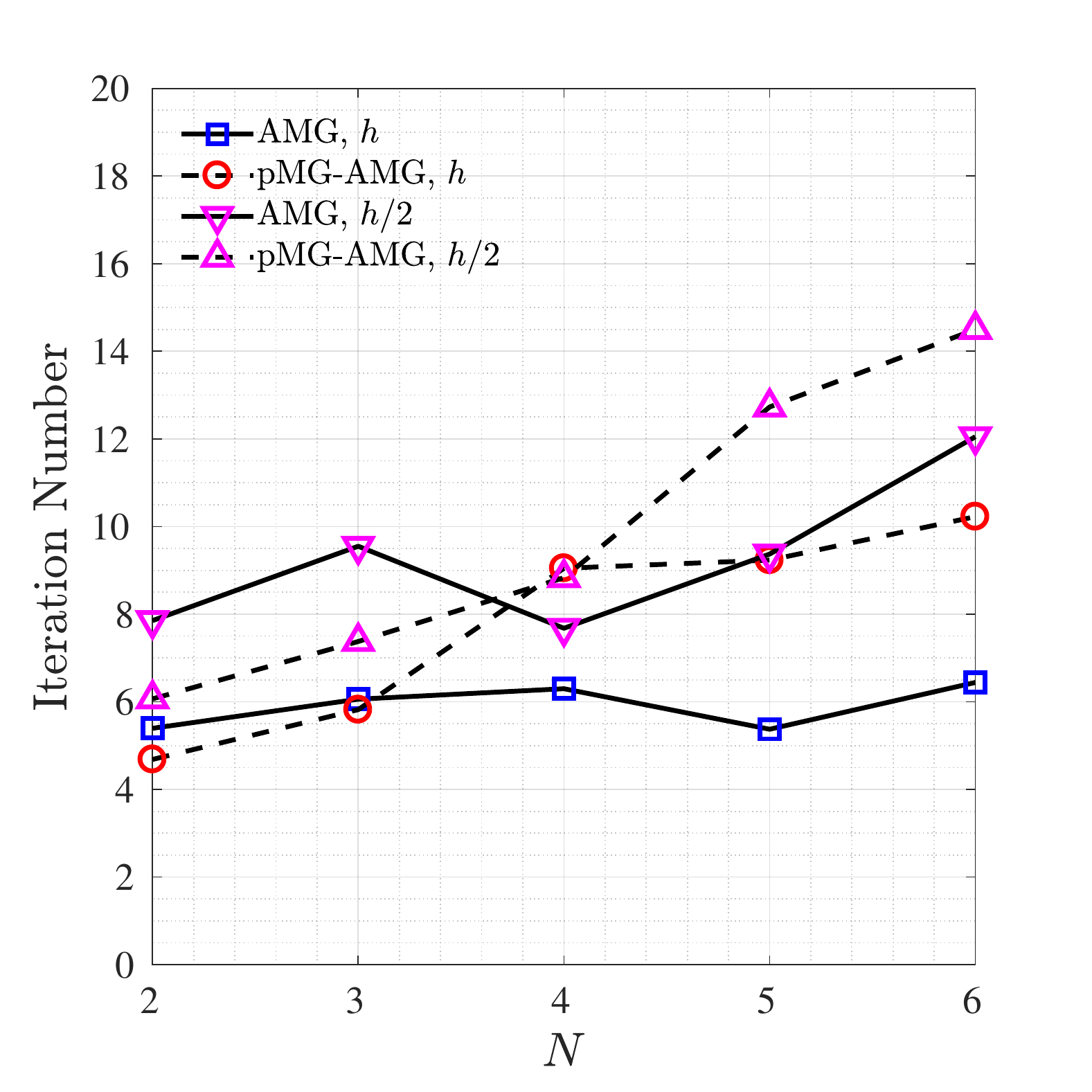}
   \caption{}
   \label{fig:VortexPreconditioner_1}
   \end{subfigure}
 \begin{subfigure}{0.45\textwidth}
  \includegraphics[width=0.99\textwidth]{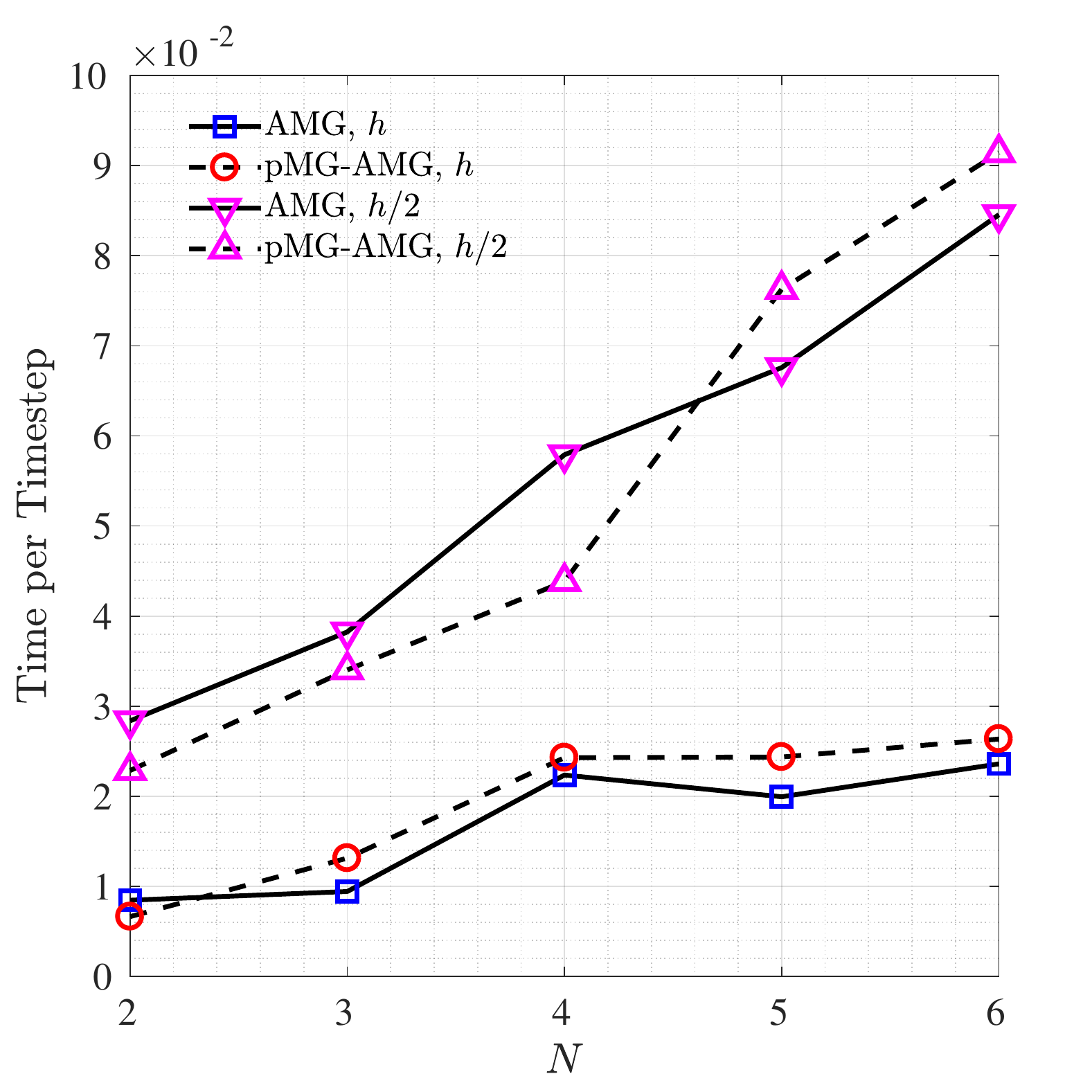}
  \caption{}
  \label{fig:VortexPreconditioner_2}
 \end{subfigure}
 \begin{subfigure}{0.45\textwidth}
  \includegraphics[width=0.99\textwidth]{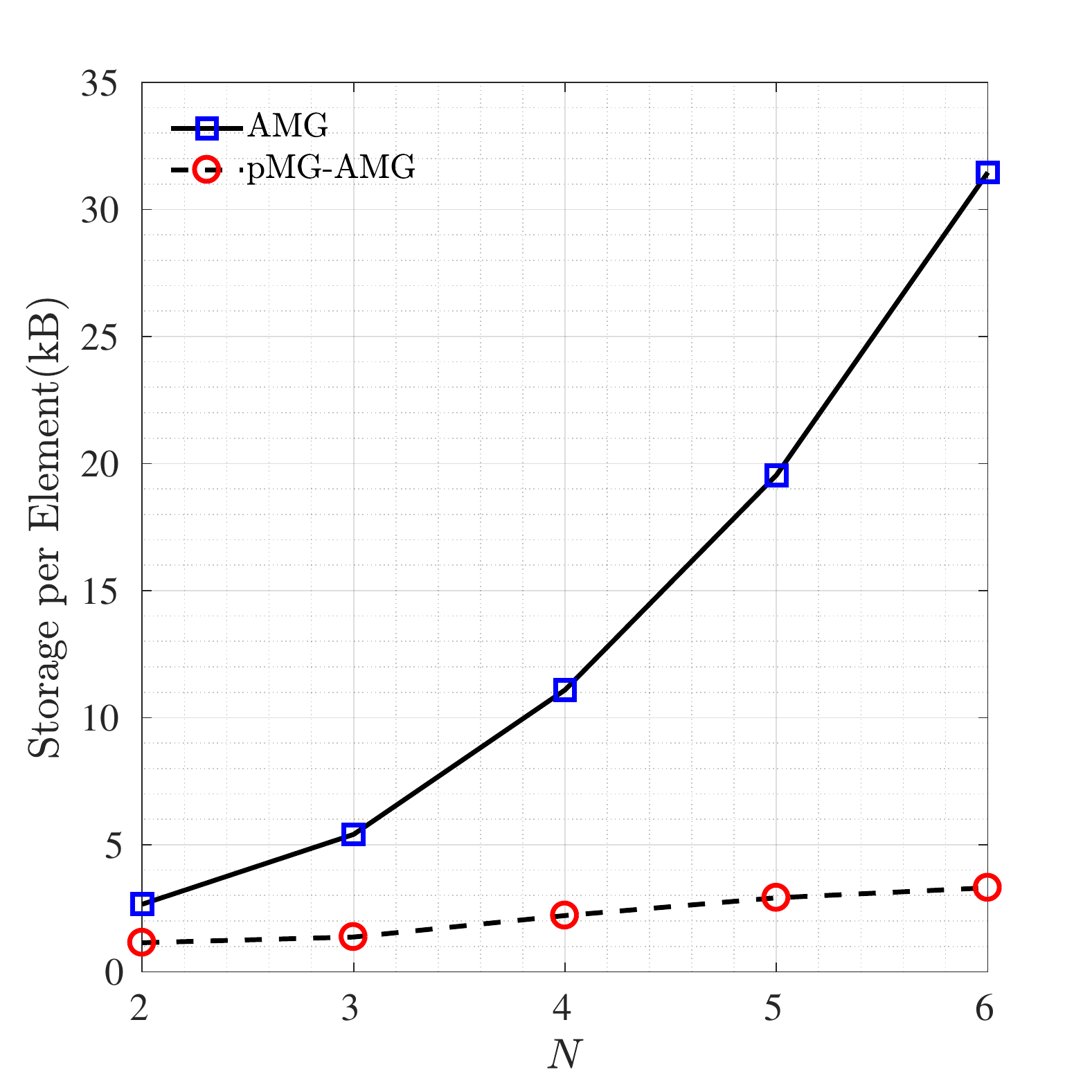}
  \caption{}
  \label{fig:VortexPreconditioner_3}
 \end{subfigure}
\caption{Comparison of hybrid $p$MG-AMG and full AMG preconditioners for Taylor vortex test problem using successively refined triangular grids at different approximation orders in terms of (a) iteration numbers (b) time spent per timestep and (c) additional storage required for the preconditioner per element.}
\label{fig:VortexPreconditioner}
\end{center}
\end{figure}
In Figure \ref{fig:VortexPreconditioner}, we compare the AMG and pMG-AMG preconditioners for the solution of pressure Poisson equation on two mesh resolutions obtained with one level uniform refinement and different approximation orders for $N=2\dots6$. For higher approximations, Figure \ref{fig:VortexPreconditioner_1} shows the number of iterations for the pMG-AMG is slightly larger than for the full AMG. Figure \ref{fig:VortexPreconditioner_2}  shows that this behavior does not lead to an increase in the time spent for each solve step. In fact, both preconditioners have comparable time-to-solution per timestep. On the other hand, the memory required for the AMG preconditioner increases dramatically with $N$. Consequently, memory requirements for the AMG preconditioner can easily exceed the limited GPU memory capacity. As shown in Figure \ref{fig:VortexPreconditioner_3}, the  AMG preconditioner uses around $30$kB of memory per element while the pMG-AMG preconditioners uses only $4$kB of storage,  and grows slowly with the order of approximation.

\subsection{Flow Past a Square Cylinder}
The relative importance of each solve step in the splitting scheme and the effect of subcycling  are  examined by solving the vortex shedding behind a square cylinder at $Re=100$. We solve the problem on a rectangular domain of size  $[-16, 25] \times [-22, 22]$ discretized with $K=2300$ unstructured triangular elements. The mesh resolution is increased near the cylinder to resolve large gradients.

\begin{figure}[htb!]
\begin{center}
   \includegraphics[width=0.95\textwidth]{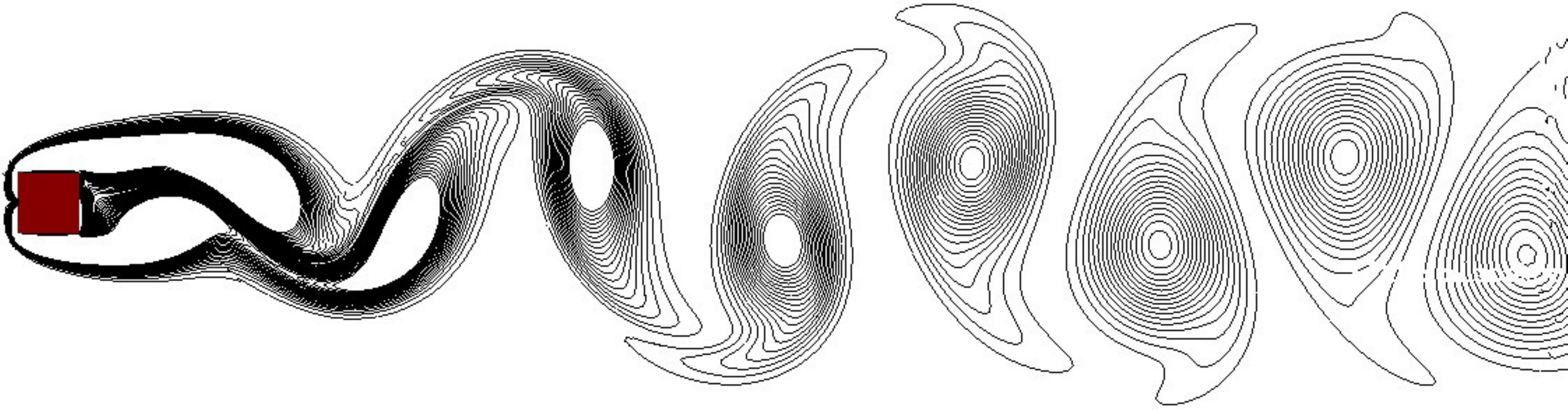}
\caption{Vortex structure in flow around square cylinder problem for $Re=100$ and approximation order, $N=5$ at time, $t=130$. Contours are from $-1$ to $1$ with increment of $0.05$.}
\label{fig:VortexShedding2D}
\end{center}
\end{figure}

The domain boundaries are inflow at the left, upper, and lower walls, outflow at the right wall, and zero Dirichlet on the square cylinder. We use zero initial conditions and unit normal velocity at inflow boundaries. Figure \ref{fig:VortexShedding2D} shows the vorticity contours of the flow at non-dimensional time $t=130$ and illustrates the instantaneous von-Karman vortex shedding profile behind the cylinder.
In order to compare our results to the available results in the
literature, we compute the Strouhal number given by $St =f D/U $,
where $f$  is the frequency of the vortex shedding, $D$ is the
characteristic length taken as the cylinder edge and $U$ is the unit characteristic velocity in this problem. We find that $St=.145$ which agrees well with the tabulated results in \citep{shahbazi_high-order_2007} and \citep{darekar_flow_2001}.

\begin{figure}[htb!]
\begin{center}
 \begin{subfigure}{0.45\textwidth}
  \includegraphics[width=0.99\textwidth]{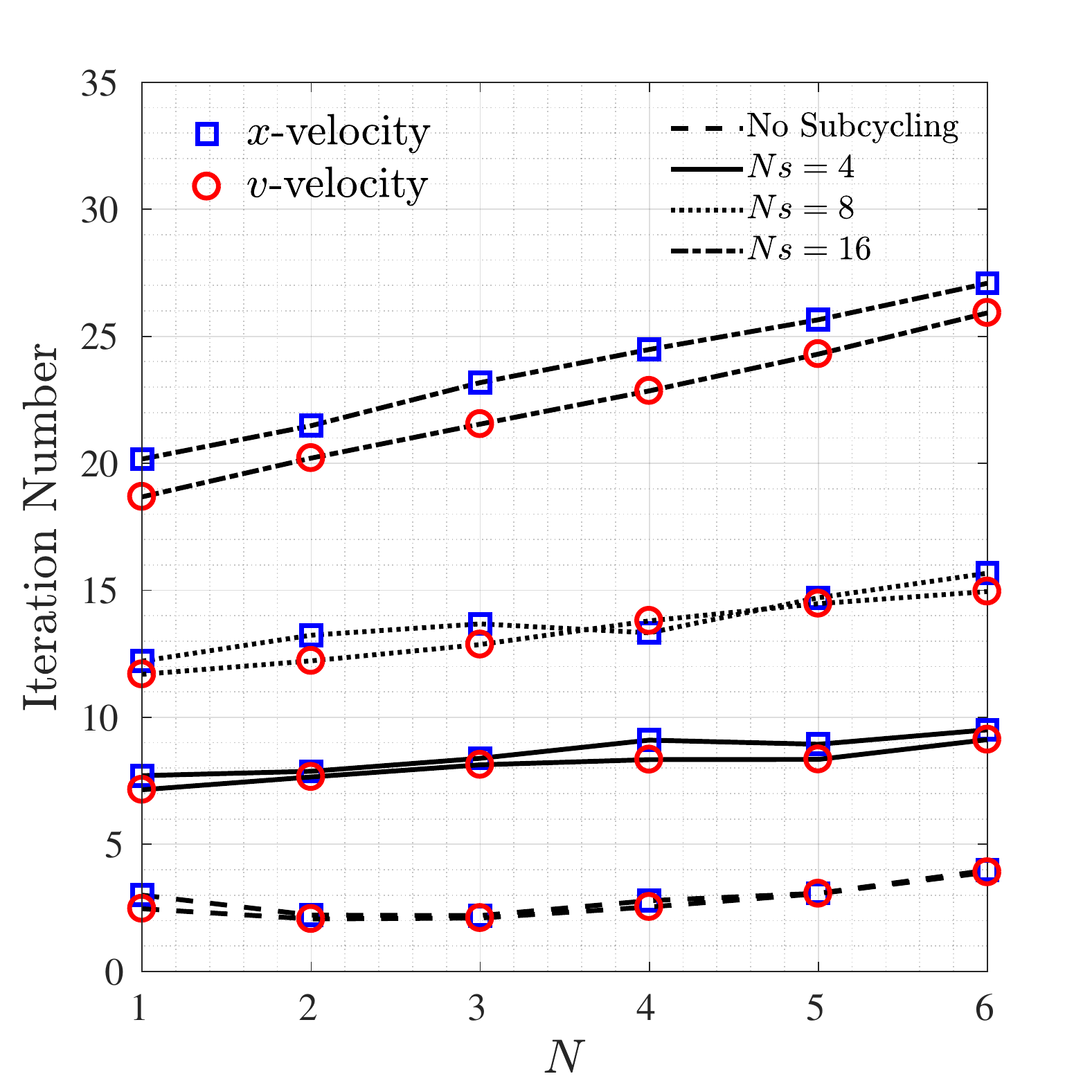}
  \caption{}
  \label{fig:SquareCylinder2DSubcycle_1}
 \end{subfigure}
 \begin{subfigure}{0.45\textwidth}
  \includegraphics[width=0.99\textwidth]{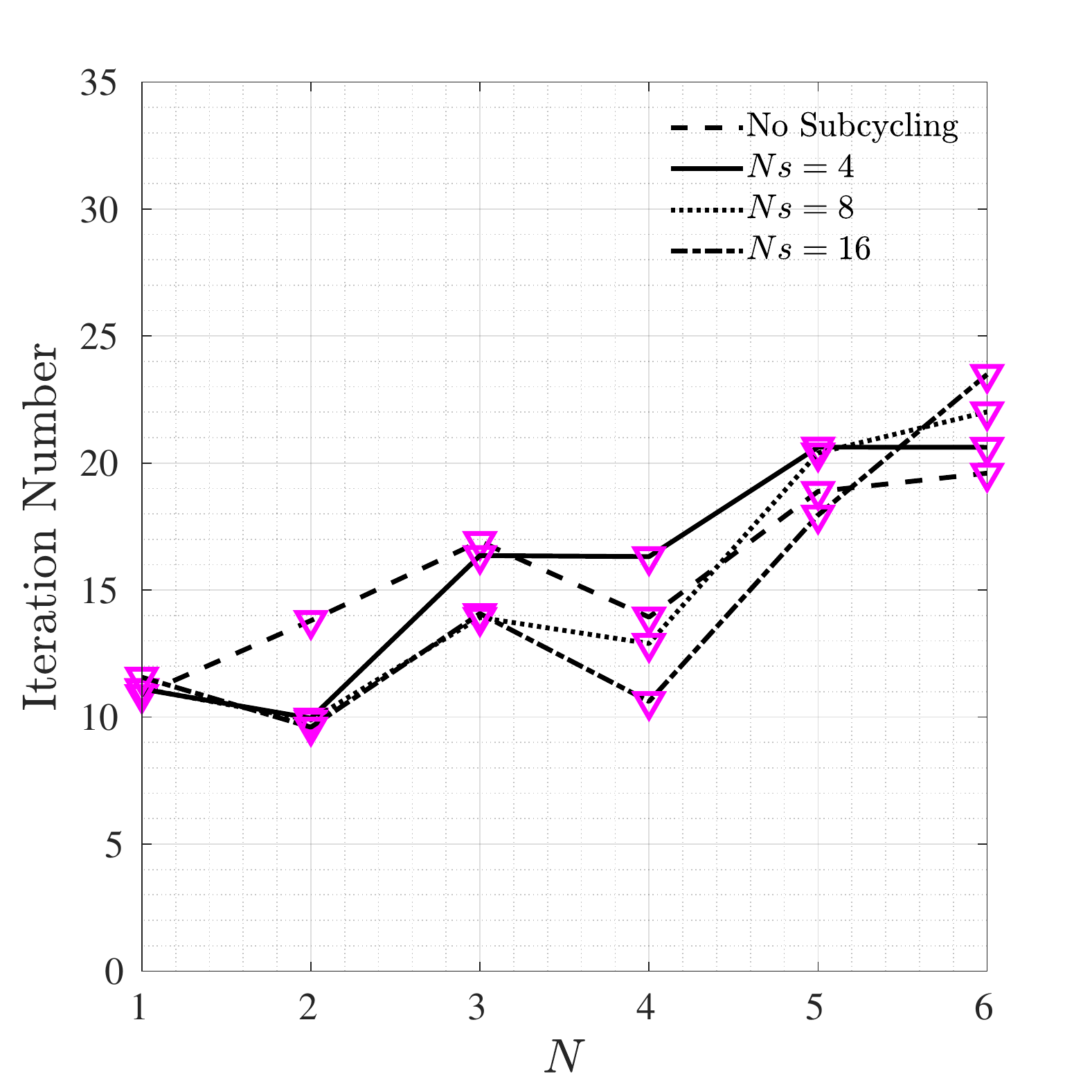}
  \caption{}
  \label{fig:SquareCylinder2DSubcycle_2}
 \end{subfigure}
 \begin{subfigure}{0.45\textwidth}
  \includegraphics[width=0.99\textwidth]{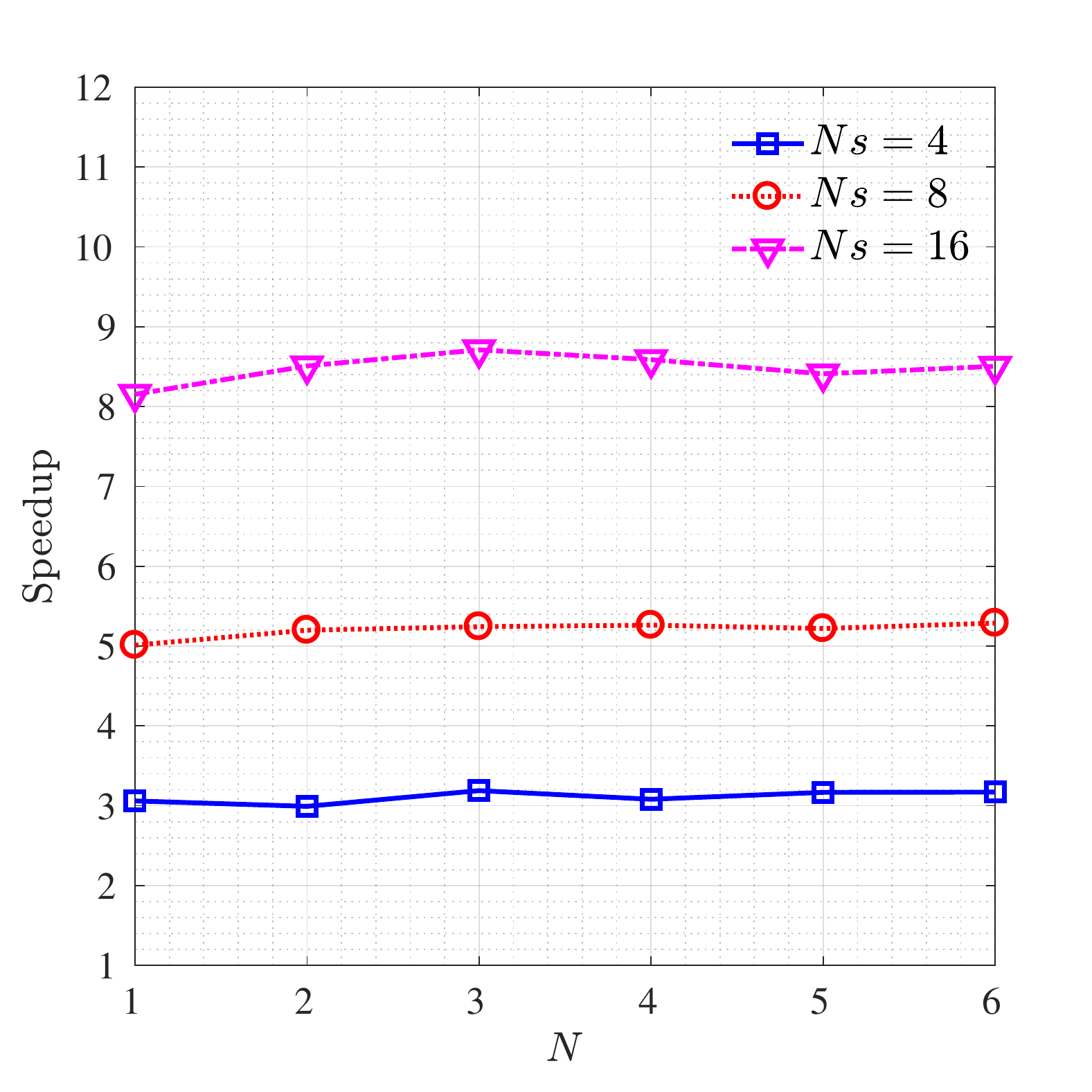}
  \caption{}
  \label{fig:SquareCylinder2DSubcycle_3}
 \end{subfigure}
\caption{Computational impact of subcycling on the individual linear solve steps and overall speedup for varying number of substeps in flow past square cylinder test case. (a) Iteration numbers in screened Poisson velocity solves, (b) Iteration numbers for pressure Poisson solve (c) Speedups using subcycling for $N_s =4$, $N_s =8$  and $N_s =16$. }
\label{fig:SquareCylinder2DSubcycle}
\end{center}
\end{figure}
Figure \ref{fig:SquareCylinder2DSubcycle} demonstrates how semi-Lagrangian subcycling affects the linear system solvers in steps \eqref{eq:INS_TD_3_2} and \eqref{eq:INS_TD_3_3}. The iteration counts required in each velocity solve are shown in Figure \ref{fig:SquareCylinder2DSubcycle_1}. We see in this figure that the iterations required increases with the number of substeps due to the larger timestep sizes making the screened Poisson operator less dominated by the mass matrix and the block-Jacobi preconditioner becoming less effective. It is important to note, however, that although iteration counts in the velocity solves are considerably higher when using subcycling, as we show below the relative time of velocity solve remains small compared with the pressure solve. Therefore the increased iteration counts do not result in an overall increase in the run times. 

On the other hand, we see in Figure \ref{fig:SquareCylinder2DSubcycle_2} that subcycling does not have an impact on the pressure solver performance. Finally, \ref{fig:SquareCylinder2DSubcycle_3} shows the achieved speedups for $N_s = 4,8,16$ and $N=1\ldots 6$. The speedups are less than the timestep size gain because of the extra computational effort required for subcycling advection step. Subcycling gives roughly $3,5$ and $8$ fold speedups for $N_s=4,8$ and $N_s=16$, respectively.

\begin{figure}[htb!]
\begin{center}
 \begin{subfigure}{0.6\textwidth}
   \includegraphics[width=0.99\textwidth]{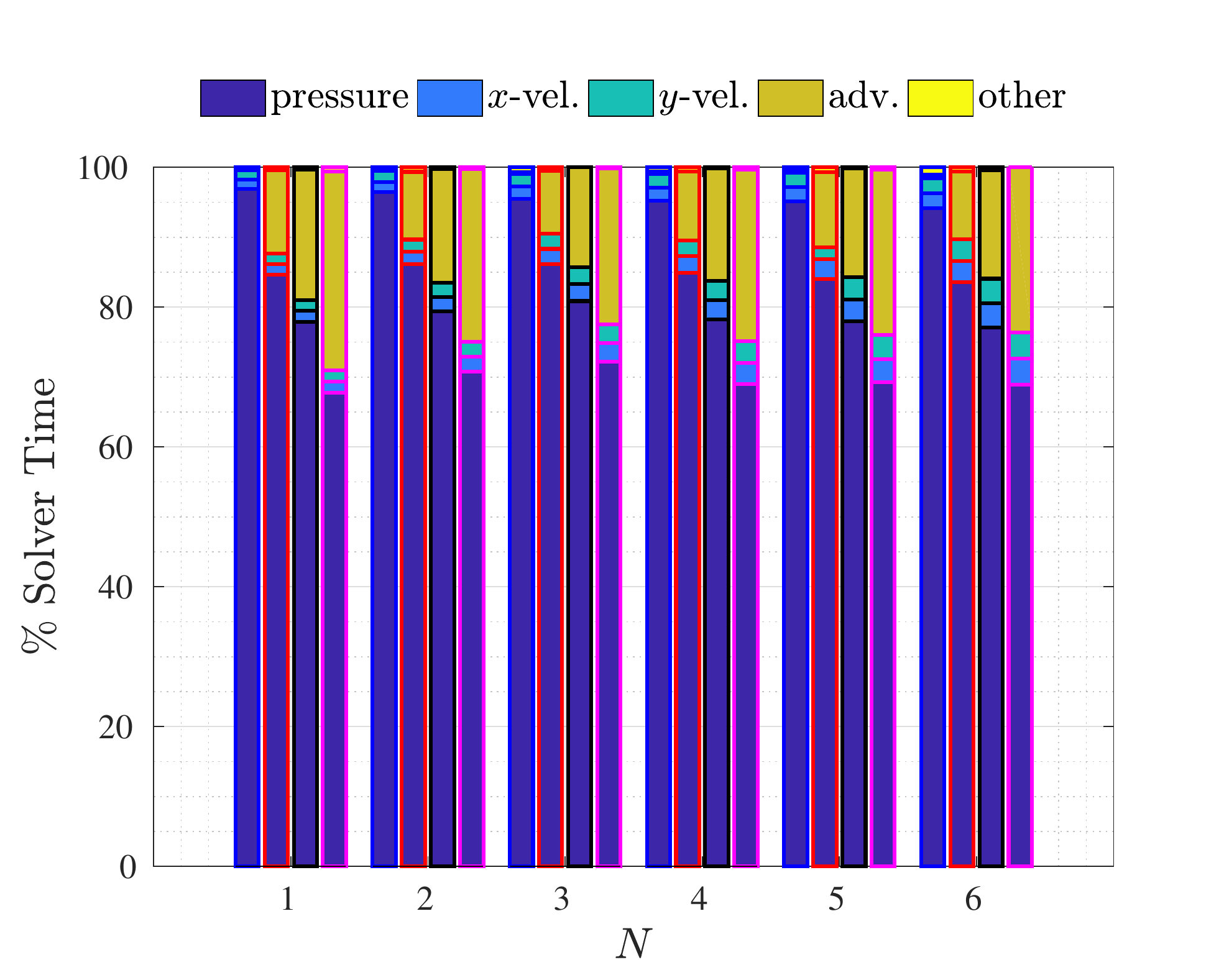}
   \caption{}
   \label{fig:SquareCyclinder2DTimings_1}
 \end{subfigure}
  \begin{subfigure}{0.6\textwidth}
   \includegraphics[width=0.99\textwidth]{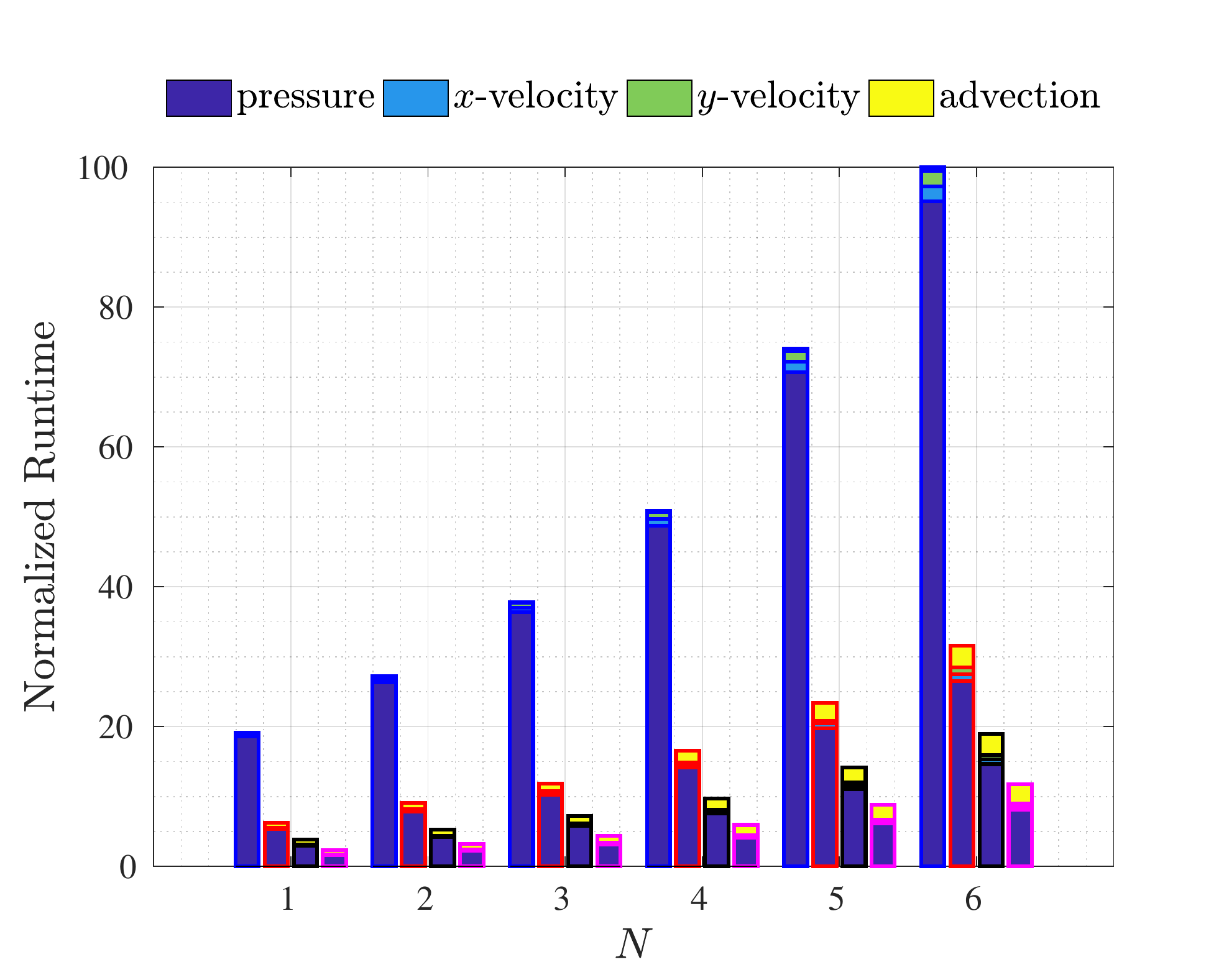}
   \caption{}
   \label{fig:SquareCyclinder2DTimings_2}
 \end{subfigure}
\caption{Timing for different treatment of  advection step (a) relative solver times (b) relative kernel times.  Each column from left to right show no-subcycling and subcycling with $N_s=4$, $N_s=8$, $N_s=16$, respectively.}
\label{fig:SquareCyclinder2DTimings}
\end{center}
\end{figure}
Figure \ref{fig:SquareCyclinder2DTimings} illustrates the percentage of time spent in each solve step, and the breakdown of normalized run times, for various numbers of subcycling steps, $N_s=4,8,16$, for  orders $N=1\ldots6$. Without subcycling, the pressure solve step takes almost all of the solution time and the overall time spent per timestep increases with the approximation order. The use of subcycling shifts the computational load away from the pressure solve to the advection steps as much more work is done in time stepping the advective terms. The resulting percentage of the time taken by the advection steps in each time step therefore becomes more significant. 

In terms of overall run times, the time taken to perform each time step of the solver decreases significantly with the use of subcycling. This is an attractive property but it requires us to give particular attention to the parallel performance of the advection kernels in optimizing the overall performance of the solver. We discuss implementation details and optimization of each of the most time consuming kernels in the next section.

\section{GPU Implementation}
The results in the last section indicate that the semi-Lagrangian subcycling method shifts computational load in each time step away from solving the Poisson problem for pressure and towards the advection stage. When considering the GPU optimization of the resulting algorithm, we have several stages and kernels which must be given specific attention. 

To test and optimize the GPU implementations of the INS solver described above we have implemented the solver using C++ together with the OCCA API and OKL kernel language \citep{medina_occa:_2014} together with MPI for distributed multi-GPU/CPU platforms. OCCA is an abstracted programming model designed to encapsulate native languages for parallel devices such as CUDA, OpenCL, Pthreads, and OpenMP. Therefore, OCCA allows customized implementations of algorithms for several computing devices with a single code and offers flexibility in choosing hardware architectures and programming model at run-time. 

For all the results presented in this section, we have compiled the source code using the GNU GCC $5.2.0$ compiler and the Nvidia CUDA V$8.0.61$ NVCC compiler. The performance tests were run using Nvidia Tesla P100 GPUs whose technical specifications are stated to be $549$ GB/s of theoretical bandwidth, $12$ GB of memory, and $4670$ GFLOPS/s of peak double-precision performance. Each GPU is running on a machine equipped with an Intel Xeon E5-$2680$v4 processor with $2.40$ GHz base frequency and $14$ cores. All the computations are performed in double precision on a fixed unstructured triangular grid with approximately $K=10,000$ elements. 

In each stage of time stepping in the INS solver we focus on the performance of the most computationally demanding kernels. In the subcycling advective stage we focus on the nonlinear volume and surface contributions of the convective term $\mathbf{u}\cdot\nabla\mathbf{u}$. Within the elliptic solve stages of each time step, in which we solve a linear system of the form $\mathbf{A}u = b$, we focus on optimizing the application of the elliptic operator $\mathbf{A}$. For the SIPDG method consists of a local gradient kernel and a kernel which computes $\mathbf{A}u$ using $u$ and $\nabla u$. 

In each section below, we give the mathematical formulation of the operators under consideration, a base pseudo-code which we implement directly in the OKL kernel language to serve as a reference implementation, and the details of successive optimizations performed to obtain better performance. We compare the GLFOPS/s achieved by each kernel version to an empirically determined roofline performance model which we detail below. 


\subsection{Empirical Roofline Model}
We evaluate the performance of our kernels by recording the run time and the number of floating point operations performed per second. Since the reported theoretical peak performance on the GPU can not be realized for most applications, we use an empirical performance model to asses the performance of our kernels. The model gives us a more realistic upper bound in terms of a maximum number of floating point operations per second that a given kernel can achieve. 

To utilize the fine-grain parallelism of the GPU we associate each thread with a single node in an element as done in \citep{klockner_nodal_2009}. This strategy has shown strong performance in previous works \citep{modave_gpu_2016}. We note, however, there exists some alternative approaches such as one thread to one element approach \citep{fuhry_discontinuous_2014}.

We consider a model to estimate the parallel performance of this strategy. Since the computational work is distributed to the individual threads on the GPU, the model is strongly based on an assumption that global data transfers and shared memory transactions limit the performance. Even if a kernel requires no floating point operations or performs only operations that are perfectly overlapped with the data movement, the runtime of this kernel cannot be shorter than the time needed to transfer the required data. 

Therefore, we consider the cost of data movement to be the most important performance limiting factor. Let us consider a kernel that loads $D_{in}$ bytes of data and stores $D_{out}$ bytes of data. We measure the time needed to transfer $ \left(D_{in} + D_{out}\right)/2$ bytes from one location in device memory to a different location. Note that we divide by $2$ due to two-way memory bus. Next, we compute a bandwidth estimate of the global memory throughput, $B_{g}$ based on the time estimate. Device to device copy bound for a kernel is determined using the formula 
\[
\frac{B_g\cdot W}{D_{in}+D_{out}},
\]
where $W$ is the work performed by the kernel, measured in GFLOPS.

We also consider the shared memory bandwidth as a supporting measure. Indeed, Volkov \citep{Volkov2008} showed that excessive shared memory read and write transactions can limit overall performance. The memory bandwidth of shared memory is estimated using the formula
\[
B_{sh}= \mathrm{\# SMs} \times \mathrm{\#ALUs} \times  \mathrm{word}\;\mathrm{length} \times \mathrm{clock}\;\mathrm{speed}\;\mathrm{in}\;\mathrm{GHz}.
\]
For the Nvidia Tesla P100 we obtain the bandwidth $B_{sh} = 7.882$ TB/s. Similar to the device to device copy bound, a shared memory performance bound can be estimated using 
\[
\frac{B_{sh}\cdot W}{S_{in}+S_{out}}, 
\]
where $S_{in}$ and $S_{out}$ are the number of bytes read and written to and from shared memory per threadblock, respectively. All the kernels considered in this section perform $4$ flops for each shared memory byte written or read. This leads to an upper bound of roughly $2$ TFLOPS/s of achievable double precision peak performance. Finally, we construct a full roofline performance model by considering the minimum of shared memory bound and device to device copy bound. 
\subsection{Elliptic Operator Kernels}
In stages \eqref{eq:INS_TD_3_2} and \eqref{eq:INS_TD_3_3} of each time step in temporal splitting scheme described above we must solve a linear elliptic system. Specifically, a screened Poisson equation for each component of the velocity field and a Poisson equation for the pressure. Optimizing solution methods of each of these systems is a difficult task, especially when considering the variety of preconditioning strategies available. In this section, we assume that the dominant cost of these linear systems is the evaluation of the elliptic operator itself. This assumption is usually well founded as iterative solution methods require several outer iterations and preconditioning methods such as multigrid require many elliptic operations at each grid level for smoothing actions.  

We detailed above the SIPDG discrete operator $L^e$ for the high-order approximation of the Laplacian operator. Here, we consider a more general operator $\mathbf{A}^e$ which approximates the screen Poisson operator on the element $\EN^e$, i.e. $\mathbf{A}^e$ approximates the action of $-\Delta + \lambda$. From the definition of $L^e$ in \eqref{eq:INS_SD_5} we can write the definition of the action of $\mathbf{A}^e$ on the polynomial $u \in V^e_N$ as satisfying
\begin{align}
(v, \mathbf{A}^e u)_{\EN^e} &= (v, -L^e u)_{\EN^e} + \lambda(v,u)_{\EN^e}, \label{eq:ellipticOp1}\\
&= (\nabla v, \nabla u)_{\EN^e} - (v,\mathbf{n}\cdot\avg{\nabla u})_{\partial \EN^e}\nonumber \\ &\qquad\qquad + \frac{1}{2}(\mathbf{n}\cdot\nabla v, \jump{u})_{\partial\EN^e} - (v,\tau\jump{u})_{\partial \EN^e} + \lambda(v,u)_{\EN^e}, \nonumber
\end{align}
for all $v \in V^e_N$. 

Next, in order to write the action of $\mathbf{A}^e$ as a linear matrix operator on the degrees of freedom of $u$ we introduce the elemental mass $\mathcal{M}^e$, surface mass $\mathcal{M}^{ef}$, and stiffness operators $\mathcal{S}^e_x$ and $\mathcal{S}^e_x$  which are defined as follows 
\begin{align}
    \mathcal{M}^e_{ij} = \left(l^e_i,l^e_j\right)_{\EN^e},  \label{eq:elMass}\quad
    \mathcal{M}^{ef}_{ij} = \left(l^e_i,l^e_j\right)_{\dE^{ef}}, 
\end{align}
\begin{align}
    (\mathcal{S}^e_x)_{ij} = \left(l_i^e,\frac{\partial l_j^e}{\partial x}\right)_{\EN^e},\label{eq:elStiff} \quad (\mathcal{S}^e_y)_{ij} = \left(l_i^e,\frac{\partial l_j^e}{\partial y}\right)_{\EN^e}.
\end{align}
Next, we define the elemental gradient operator $\boldsymbol{\mathcal{D}}^e = [\mathcal{D}^e_x, \mathcal{D}^e_y]^T$, as well as the lifting operators $\mathcal{L}^{ef}$, via
\begin{align}
    \mathcal{D}^e_x = (\mathcal{M}^e)^{-1}\mathcal{S}_x,\quad
    \mathcal{D}^e_y = (\mathcal{M}^e)^{-1}\mathcal{S}_y, \quad
    \mathcal{L}^{ef} = (\mathcal{M}^e)^{-1}\mathcal{M}^{ef}. \label{eq:elLift}    
\end{align}
Finally, for ease of notation we introduce the concatenation of the lift operators along each face, i.e. $\mathcal{L}^{e} = [\mathcal{L}^{e0}, \mathcal{L}^{e1}, \mathcal{L}^{e2}]$.

Returning to the elliptic operator $\mathbf{A}^e$ in \eqref{eq:ellipticOp1}, to improve the performance we aim to avoid using transpose versions of the operators defined above. We also aim to avoid performing excessive matrix-vector products. To this end, we rewrite this operator to group common operations as much as possible. To begin, we integrate the first volume integral in the expression above to obtain 
\begin{multline}\label{eq:ellipticOp2}
(v, \mathbf{A}^e u)_{\EN^e} = -(v, \Delta u)_{\EN^e} - \frac{1}{2}(v,\mathbf{n}\cdot\jump{\nabla u})_{\partial \EN^e} \\ + \frac{1}{2}(\mathbf{n}\cdot\nabla v, \jump{u})_{\partial\EN^e} - (v,\tau\jump{u})_{\partial \EN^e} + \lambda(v,u)_{\EN^e}.   
\end{multline}
Next, we note that from the from the definition of the lift operators $\mathcal{L}^{ef}$ in \eqref{eq:elLift} we can write 
\begin{align*}
    (\mathbf{n}\cdot\nabla v, \jump{u})_{\partial\EN^e} &= (\mathbf{n}\cdot\nabla v, \mathcal{L}^{e}\jump{u})_{\EN^{e}}, \\
    &= -(v, \mathbf{n}\cdot\nabla\mathcal{L}^{e}\jump{u})_{\EN^{e}} + (v, (\mathcal{L}^{e}\jump{u})^-)_{\partial\EN^{e}}.
\end{align*}
Here we applied integration by parts to obtain the last line, recalling that the $-$ superscript denotes the interior trace. Using this expansion in \eqref{eq:ellipticOp1} we obtain 
\begin{multline*}
(v, \mathbf{A}^e u)_{\EN^e} = -(v, \Delta u)_{\EN^e} - \frac{1}{2}(v,\mathbf{n}\cdot\jump{\nabla u})_{\partial \EN^e} - \frac{1}{2}(v, \mathbf{n}\cdot\nabla\mathcal{L}^{e}\jump{u})_{\EN^{e}}\\ + \frac{1}{2}(v, (\mathcal{L}^{e}\jump{u})^-)_{\partial\EN^{e}} - (v,\tau\jump{u})_{\partial \EN^e} + \lambda(v,u)_{\EN^e}.   
\end{multline*}
Finally, taking $v$ to be each of the basis polynomials $l^e_n$, $n=1,\ldots,N_p$, we can use the elemental operators defined in \eqref{eq:elMass}-\eqref{eq:elLift} in order to write the action of operator $\mathbf{A}^e$ on the polynomial $u$ as 
\begin{align}
    \mathbf{A}^e u &= -\mathcal{M}^e\boldsymbol{\mathcal{D}}^e\cdot\boldsymbol{\mathcal{D}}^e u -\frac{1}{2}\mathcal{M}^e\mathcal{L}^{e}\mathbf{n}\cdot\jump{\boldsymbol{\mathcal{D}}^e u}    -\frac{1}{2}\mathcal{M}^e\mathbf{n}\cdot\boldsymbol{\mathcal{D}}^e\mathcal{L}^{e}\jump{u} \nonumber\\
    &\hspace{5cm}
    +\frac{1}{2}\mathcal{M}^e\mathcal{L}^{e}(\mathcal{L}^{e}\jump{u})^-
    -\tau\mathcal{M}^e\mathcal{L}^{e}\jump{u} + \lambda\mathcal{M}^e u, \nonumber\\
    &= \mathcal{M}^e\left(-\boldsymbol{\mathcal{D}}^e\cdot\left[\boldsymbol{\mathcal{D}}^e u +\frac{1}{2}\mathbf{n}\mathcal{L}^e\jump{u}\right] - \frac{1}{2}\mathcal{L}^{e}\left[ \mathbf{n}\cdot\jump{\boldsymbol{\mathcal{D}}^e u} + 2\tau\jump{u} - (\mathcal{L}^{e}\jump{u})^- \right] \right). \label{eq:ellipticOp3} 
\end{align}
We use expression \eqref{eq:ellipticOp3} as a basis for implementing the action of the elliptic operator $\mathbf{A}$. 

In order to obtain a more unified expression for the action of $\mathbf{A}^e$ between separate elements we introduce a mapping from each element $\EN^e$ to a reference element $\hat{\EN}$, on which we make use of reference operators. We take the reference element $\hat{\EN}$ to be the bi-unit triangle 
\[
\hat{\EN} = \left\{ -1\leq r,s,r+s\leq 1 \right\},
\]
and introduce the affine mapping $\Phi^e$ which maps $\EN^e$ to a reference triangle $\hat{\EN}$, i.e.
\begin{equation}
    \label{eq:operators1}
    \left(x,y\right) =\Phi^e\left(r,s\right), \quad \left(x,y\right)\in \EN^e,\ \left(r,s\right)\in\hat{\EN}
\end{equation}
We denote the Jacobian of this mapping as
\begin{equation}
    \label{eq:operators2}
    G^e = \begin{bmatrix}
    r_x & s_x \\
    r_y & s_y 
    \end{bmatrix},
\end{equation}
and denote determinant of the Jacobian as $J^e = \det G^e$. We also define the surface scaling factor $J^{ef}$ which is defined as the determinant of the Jacobian $G^e$ restricted to the face $\partial EN^{ef}$.  

Finally, mapping each of the elemental operators defined in \eqref{eq:elMass}-\eqref{eq:elLift} to the reference element $\hat{\EN}$ we can write each of the elemental operators in terms of their reference versions and the geometric factors $G^e$, $J^e$, and $J^{ef}$ as follows
\begin{align}
    \mathcal{M}^e = J^e \mathcal{M}, \quad 
    \boldsymbol{\mathcal{D}}^e = G^e \boldsymbol{\mathcal{D}},\quad \label{eq:elementOps}
    \mathcal{L}^{ef} = \frac{J^{ef}}{J^e}\mathcal{L}^f.
\end{align}
Here $\mathcal{M}$, $\boldsymbol{\mathcal{D}} = [\mathcal{D}_r,\mathcal{D}_s]^T$, and $\mathcal{L}^f$ are the mass, derivative, and lifting operators defined on the reference element $\hat{\EN}$. Therefore, we can write the elliptic operator $\eqref{eq:ellipticOp3}$ on each element using only these reference operators and the geometric data $G^e$, $J^e$, and $J^{ef}$. 

\subsubsection{Local Gradient Kernel}
To implement the elliptic operator on the GPU we first note that since the we require the positive and negative traces of the local derivative term $\boldsymbol{\mathcal{D}}^e u$ we must first compute and store it in global device memory so each element's neighbour data is visible. To perform this operation we first implement a local gradient kernel which inputs a field $u$ and outputs the local gradient $\boldsymbol{\mathcal{D}}^e u$.  We give the pseudo-code of this kernel in Algorithm~\ref{alg:AxG}. Since the size of the matrix-vector products in this kernel are $N_p\times N_p$ we launch this kernel using $N_p$ threads per block. 

\begin{algorithm}[t]
  \caption{Local Gradient Kernel}
  \label{alg:AxG}
\begin{boxedminipage}{1.0\textwidth}
    \begin{algorithmic}[1]
    \STATE {\bf Input:} 
    \\(1) $u$, size  $K \times N_p$.
    \\(2) Derivative  matrices $\boldsymbol{\mathcal{D}} = [\mathcal{D}_r, \mathcal{D}_s]$, size $2\times \left(N_p\times N_p\right)$.
    \\(3) Geometric factors $G$, size $ 4\times K$.
    \STATE {{\bf Output: } \\$\nabla u = [u_x, u_y]$, size  $2 \times \left(K\times N_p\right)$}.
     \FOR {$e\in\left\{1,2, \ldots K \right\}$}
     \FOR{$i\in\left\{1,2, \ldots N_p\right\}$}
     \STATE $u_{r;i}=\sum_{j=1}^{N_p}\mathcal{D}_{r;ij} u_{j}^{e}$ \COMMENT{Apply reference derivatives}
     \STATE $u_{s;i}=\sum_{j=1}^{N_p}\mathcal{D}_{s;ij} u_{j}^{e}$
     \STATE $r_{x} = G_{0}^{(e)},\; s_{x} = G_{1}^{(e)} $ \;$r_{y} = G_{2}^{(e)},\; s_{y} = G_{3}^{(e)} $
     \STATE $u_{x;i}^{e} = r_{x}\phi_{r;i} + s_{x}\phi_{s;i}$  \COMMENT{Apply geometric factors}  
     \STATE $u_{y;i}^{e} = r_{y}\phi_{r;i} + s_{y}\phi_{s;i}$
    \ENDFOR
    \ENDFOR
  \end{algorithmic}
\end{boxedminipage}
\end{algorithm}

\begin{figure}[tb]
\begin{center}
  \includegraphics[width=0.9\textwidth]{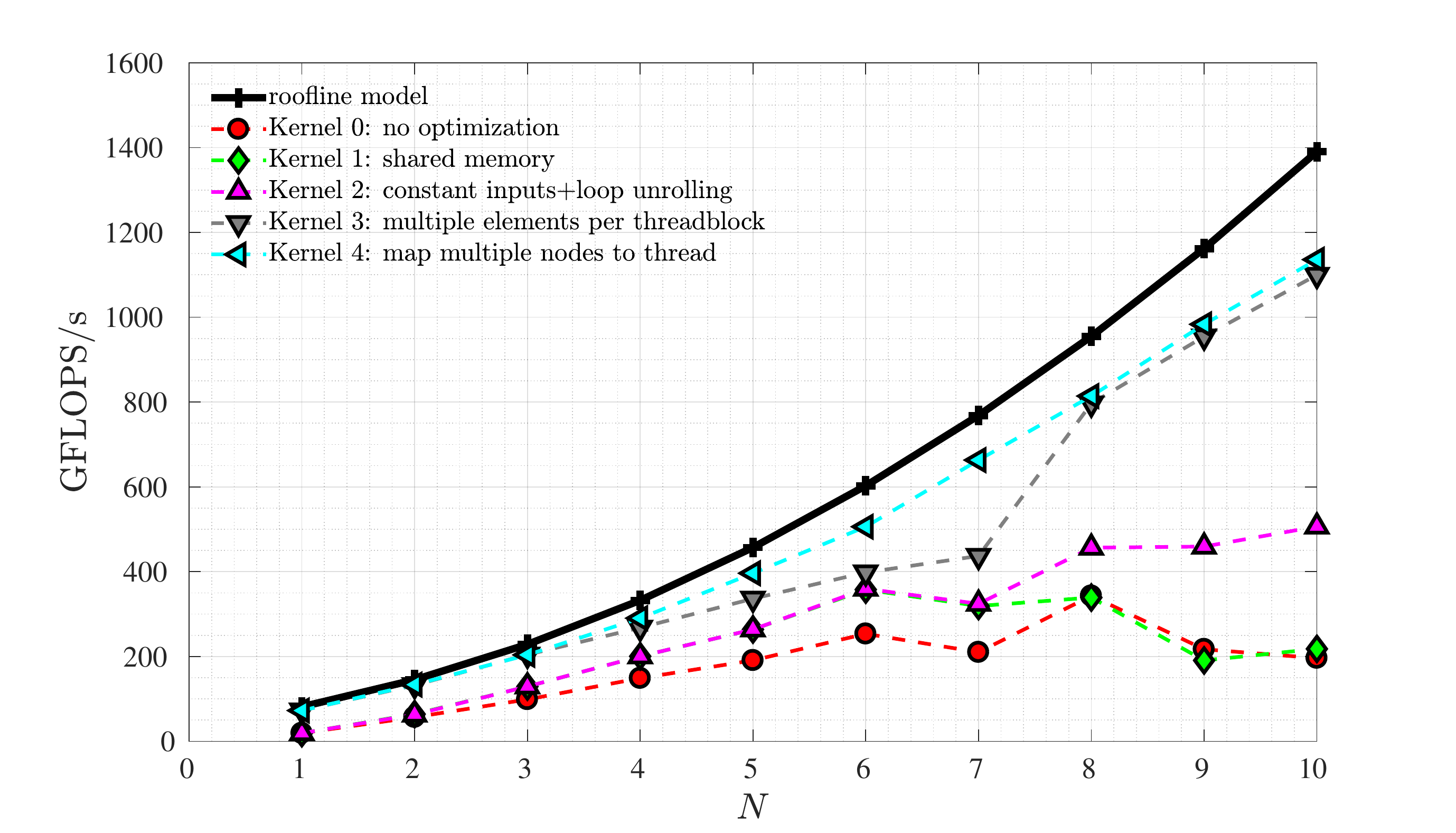}
\caption{Achieved floating point performance for the local gradient kernels compared against the empirical roofline model shown as a black line.}
\label{fig:GradientKernel}
\end{center}
\end{figure}

We show in Figure \ref{fig:GradientKernel} the GPU performance results of five kernels implementing the local gradient operation. The kernels are constructed in a sequential fashion starting with a direct implementation of Algorithm \ref{alg:SSV} and applying successive optimizations. Each kernel uses the previous kernel implementation as a starting point and applies the optimizations detailed below. 

\kernelNl{Local Gradient}{0} This kernel is a direct implementation of Algorithm~\ref{alg:AxG}. The kernel reads the $u$ field directly from global GPU memory during the matrix-vector product with the differentiation matrices. Due to these excessive global memory transactions, this kernel only reaches $200$ GFLOPS/s.

\kernelNl{Local Gradient}{1} In this kernel we add two shared memory arrays of size $N_p$ to store the $u$ field before differentiation. Using shared memory rather than repeated accesses to global memory improves the performance substantially for $N<7$. However, at higher orders the performance stalls.  

\kernelNl{Local Gradient}{2} In this kernel all the global and local variables that are not modified are labeled with $\texttt{const}$ qualifier. Also, the $\texttt{restrict}$ qualifier is added to all input arrays to indicate to the compiler that memory locations pointed to do not overlap. Furthermore, all serial loops in the differentiation actions are unrolled, increasing instruction-level parallelism. These optmizations improve the performance of the kernel for high-order approximations and the performance reaches approximately $500$ GFLOPS/s.       

\kernelNl{Local Gradient}{3} In this kernel multiple elements are processed by each threadblock to better align the computational load with the hardware architecture. Running several trials, we choose the number of elements per threadblock which optimizes performance. This optimization strategy increases the performance marginally. Achieved performance reaches $1.1$ TFLOPS/s at $N=10$ but remains below the empirical bound.

\kernelNl{Local Gradient}{4} In this kernel, each thread processes multiple nodes of an element, in addition to each threadblock processing multiple elements. That is, each time an entry of the differentiation matrices is loaded from memory it can be reused multiple times in the matrix-vector multiplication. The results of the matrix-vector products are stored in a register array. With this optimization strategy, overall performance curve of the kernel approaches the roofline curve for $N<8$. For higher order, the difference between the achieved and empirical roofline performance is approximately $10\%$.

\subsubsection{SIPDG Operator Kernel}
Once the local gradient of the field $u$ is computed and stored in global memory we use the SIPDG operator kernel to compute the action of the $\mathbf{A}$ operator on the field $u$. We give the pseudo-code of this kernel in \ref{alg:AxKernel}. As an input to this kernel we assume that an index array of negative and positive trace indices has been constructed.

\begin{algorithm}[htbp!]
  \caption{SIPDG Kernel}
  \label{alg:AxKernel}
  \begin{boxedminipage}{1.0\textwidth}
    \begin{algorithmic}[1]
    \STATE {\bf Input:} \\(1) $u$, size  $K \times N_p$.
    \\(2) $\nabla u = [u_x, u_y]$, size  $2\times (K \times N_p)$;
    \\(3) Negative trace indices $idM$, size $K \times  N_f \times N_{fp} $.
    \\(4) positive trace indices $idP$, size $K \times  N_f \times N_{fp} $.
    \\(2) Derivative  matrices $\boldsymbol{\mathcal{D}} = [\mathcal{D}_r, \mathcal{D}_s]$, size $2\times \left(N_p\times N_p\right)$.
    \\(5) Lift matrix $\mathcal{L}$, size $N_p \times N_f \times N_{fp}$. 
    \\(6) Mass Matrix $\mathcal{M}$, size $N_p\times N_p$.
    \\(7) Surface geometric factors $sG$, size $ 4\times K \times N_f$.
    \\(8) Volume geometric factors $G$, size $ 5\times K$.
    \STATE {{\bf Output: } $\mathbf{A}u$, size  $K\times N_p$}.
\FOR {$e\in\left\{1,2, \ldots K \right\}$}
  \FOR{$i\in\left\{1,2, \ldots \max\left(N_p,\left( N_f \times N_{fp}\right)\right)\right\}$}
    
     \IF[Load data and lift  jumps]{$i\leq N_p$}
      \STATE $n_{x} = sG_{0}^{e,f},\; n_{y} = sG_{1}^{e,f}  \;J^{e,f} = sG_{2}^{e,f},\;(J^{e})^{-1} = sG_{3}^{e,f}  $
      \STATE $L_{x;i}= 0.5 n_x J^{e,f} (J^{e})^{-1}\sum_{j=1}^{N_f \times N_{fp}} \mathcal{L}_{ij} \left(u(idP^{e}_j) - u(idM^{e}_j) \right)$
      \STATE $L_{y;i}=0.5 n_x J^{e,f} (J^{e})^{-1}\sum_{j=1}^{N_f \times N_{fp}} \mathcal{L}_{ij} \left(u(idP^{e}_j) - u(idM^{e}_j) \right)$
      \ENDIF
    \ENDFOR
      
    \FOR{$i\in\left\{1,2, \ldots \max\left(N_p,\left( N_f \times N_{fp}\right)\right)\right\}$}  
      \IF[Compute volume contribution]{$i\leq N_p$}
      \STATE $r_{x} = G_{0}^{(e)},\; s_{x} = G_{1}^{e} $ \;$r_{y} = G_{2}^{e},\; s_{y} = G_{3}^{e}, \; J^{e} = G_{4}^{e} $
      \STATE $Au_{i} \;\; = -\sum_{j=1}^{N_p} \mathcal{D}_{r;ij} (r_{x}(u_{x;i}^e + L_{x;i}) + r_y (u_{y;i}^e + L_{y;i}))$
      \STATE $Au_{i} \mathrel{{+}{=}} -\sum_{j=1}^{N_p} \mathcal{D}_{s;ij} (s_{x}(u_{x;i}^e + L_{x;i}) + s_y (u_{y;i}^e + L_{y;i}))$
      \ENDIF
      \IF[Compute surface contributions]{$i\leq N_f \times N_{fp}$} 
      \STATE $s_{i} \;\; = 0.5 J^{e,f}(J^{e})^{-1} n_x \left(u_x(idP^e_i) - u_x(idM^e_i)\right) $ 
      \STATE $s_{i} \mathrel{{+}{=}} 0.5 J^{e,f}(J^{e})^{-1} n_y\left(u_y(idP^e_i) - u_y(idM^e_i) \right)$ 
      \STATE $s_{i} \mathrel{{+}{=}} J^{e,f}(J^{e})^{-1}\tau \left(u(idP^e_i)- u(idM^e_i)\right)$
      \STATE $s_{i} \mathrel{{-}{=}} J^{e,f}(J^{e})^{-1}\left(n_x L_{x}(idM^e_i) + n_y L_{y}(idM^e_i) \right)$
      \ENDIF
    \ENDFOR  
     
     \FOR{$i\in\left\{1,2, \ldots \max\left(N_p,\left( N_f \times N_{fp}\right)\right)\right\}$}  
      \IF[lift surface contribution]{$i\leq N_p$}
      \STATE $Au_{i} \mathrel{{-}{=}} \sum_{j=1}^{N_f\times N_{fp}} \mathcal{L}_{ij} s_{j}$
      \ENDIF
    \ENDFOR  
      
    \FOR{$i\in\left\{1,2, \ldots \max\left(N_p,\left( N_f \times N_{fp}\right)\right)\right\}$}  
    \IF[Multiply with mass matrix]{$i\leq N_p$}
    \STATE $\mathbf{A}u_{i}^e=J^{e}\sum_{j=1}^{Np}\mathcal{M}_{ij} Au_j$
    \ENDIF
    \ENDFOR
    \ENDFOR
  \end{algorithmic}
\end{boxedminipage}
\end{algorithm}

To fully paralleize the kernel we require $N_f \times  N_{fp}$ threads for the surface flux construction, where $N_f$ is the number of faces per element and $N_{fp}$ is the number of degrees of freedom per face, and we require $N_p$ threads to paralleize the derivative and lifting operations. Therefore, we use a total of $\max\left(N_f\times N_fp, Np\right)$ threads per block with this kernel. 

\begin{figure}[tb]
\begin{center}
  \includegraphics[width=0.9\textwidth]{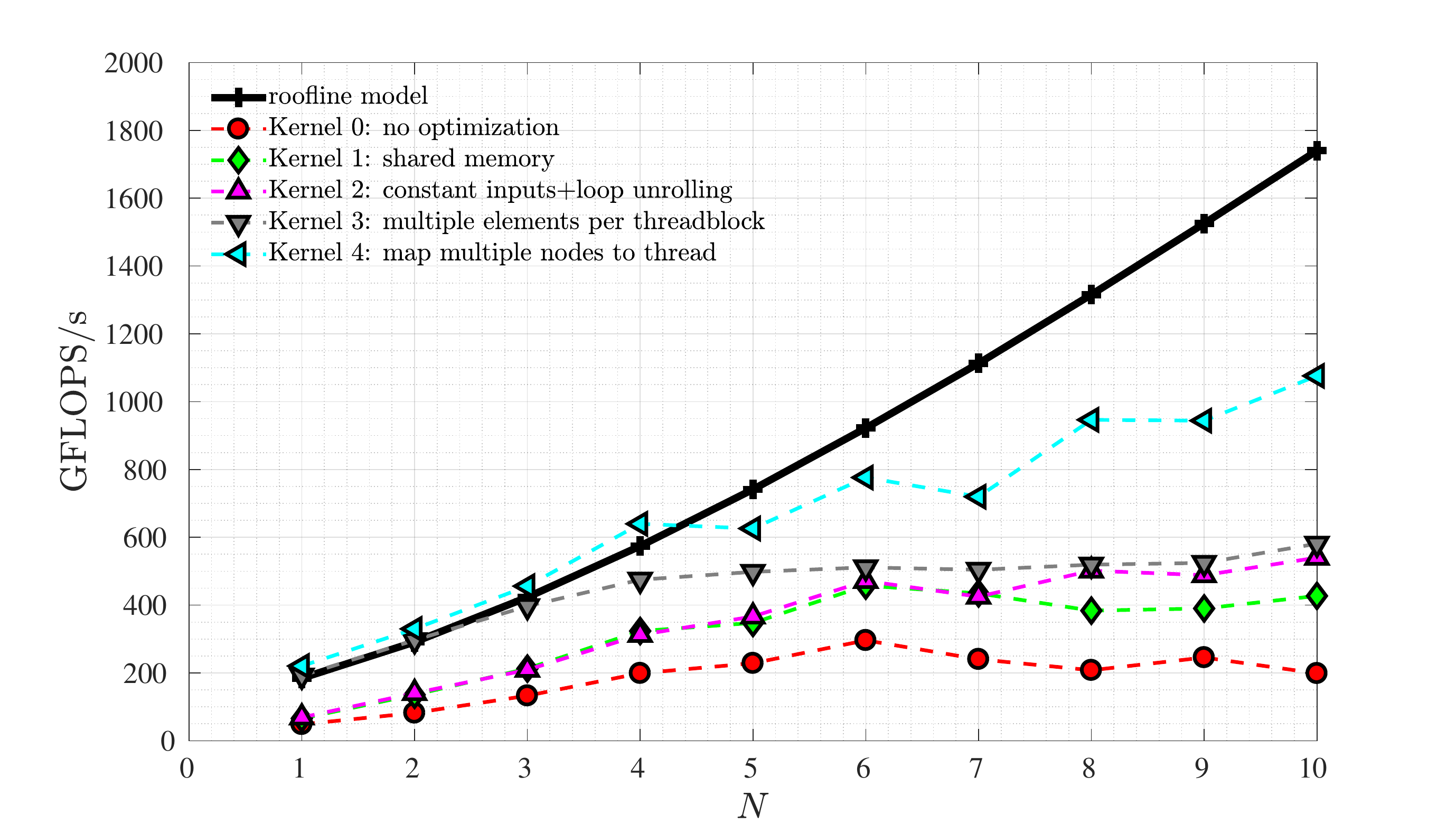}
\caption{Achieved floating point performance for the Ax kernels compared against the empirical roofline model shown as a black line. }
\label{fig:AxKernel}
\end{center}
\end{figure}

We show in Figure \ref{fig:AxKernel} the GPU performance results of five kernels implementing the SIPDG elliptic operator. As before, the kernels are constructed in a sequential fashion starting with a direct implementation of Algorithm \ref{alg:AxKernel} and applying successive optimizations. Each kernel uses the previous kernel implementation as a starting point and applies the optimizations detailed below.

\kernelNl{SIPDG}{0} This kernel is a direct implementation of Algorithm~\ref{alg:AxKernel}. In this kernel the field variable $u$ and the derivative $u_x$ and $u_y$ are loaded from global memory in lifting, volume, and surface evaluation steps. Results from matrix-vector products are stored in separate shared memory arrays. Due to the excessive global memory reads, this kernel achieves only 200 GFLOPS/s at $N>4$ which is $10\%$ of achievable performance for $N=10$. 

\kernelNl{SIPDG}{1} In this kernel we use five shared memory arrays of size $N_p$ and $N_f \times N_{fp}$ to store $u_x, u_y$ and the local and external trace values of  $u_x, u_y$, and $u$. All trace data is loaded from global memory before first lifting step, which requires a thread synchronization to ensure cache coherence. Reducing the global memory transactions increases the performance of this kernel by roughly a factor of two. 

\kernelNl{SIPDG}{2} In this kernel we add a \texttt{const} qualifier to all input and local variable which remain unmodified and add the \texttt{restrict} qualifier to all input arrays. We also unroll serial $\texttt{for}$ loops to increase instruction-level parallelism. This kernel reaches $550$ TFLOPS/s for $N>6$ but we do not see a significant improvement for lower orders.  

\kernelNl{SIPDG}{3} In this kernel multiple elements are processed by each threadblock to increase occupancy. The number of elements mapped to a threadblock is optimized for each order of approximation by running several trials. Performance of the kernel increases substantially for low orders, and the measured performance approaches the empirical roofline curve. For $N>4$ achieved performance stalls around $600$ GFLOPS/s. This behavior can be explained by excessive operator loads. The SIPDG kernel requires a mass matrix, lift operator, and local differentiation matrices with sizes  $N_p\times N_p$, $N_p \times \left(N_f\times N_fp\right)$ and $2\times\left(N_p\times N_p\right)$, respectively. For $N>4$, the data fetched by the kernel exceeds $24$KB, which is the capacity of L1 cache in an Nvidia Tesla P100 GPU. Since these operators cannot be stored in cache for $N>4$, we observe a drop in performance due to global memory cache-misses. 

\kernelNl{SIPDG}{4} In this kernel, in addition to processing multiple elements in a threadblock, multiple nodes are processed by a single thread. This strategy allows for reusing operators multiple times per load and, hence, brings considerable performance improvement. The observed performance curve approaches the roofline curve for low orders and reaches $1.1$ TFLOPS/s at $N=10$ with less observed stalling for $N>4$. The kernel still achieves only $30\%$ of the predicted achievable performance due to L1 cache misses, and nonsequential data access pattern of external trace values leading to reduced data coalescing in global reads.

\subsection{Subcycling Advection Kernels}
The subcycling method requires several evaluations of the operator $\tilde{\mathbf{N}}^e(\bar{\mathbf{U}}^e,\tilde{\mathbf{U}}^e)$ defined in \eqref{eq:INS_CUB_N}. To describe the evaluation of this operator, we first note that we must use a sufficiently high-order cubature rule to prevent aliasing errors when evaluating the integrals in \eqref{eq:INS_CUB_N}. We consider a sufficient  nodal set of $N_c$ cubature nodes with coordinates in the reference element $(r^c_i,s^c_i)$, and associated weight, $w^c_i$ for $i=1,\ldots,N_c$. We define analogous cubature node set on each face of $\hat{\EN}$ in order to integrate the surface terms with sufficiently high-order and for each face we denote these nodes as $(r^{cf}_j,s^{cf}_j)$, and associated weight, $w^{cf}_j$ for $j=1,\ldots,N^f_c$. We use these cubature nodes to define the following interpolation operators
\begin{align}
    (\mathcal{I} \phi(r,s))_i &= \phi(r_i^c, s_i^c), \label{eq:operators5} \\
    (\mathcal{I}^f \phi(r,s))_j &= \phi(r_j^{cf}, s_j^{cf}),
\end{align}
for $i=1,\ldots,N_c$ and $j=1,\ldots,N_c^f$. 

Mapping \eqref{eq:INS_CUB_N} to the reference element $\hat{\EN}$ and taking the test functions $v$ to be each of the nodal basis functions $v= l_n$ we find that we can write the operator $\tilde{\mathbf{N}}^e(\bar{\mathbf{U}}^e,\tilde{\mathbf{U}}^e)$ as
\begin{multline*}
(J^e \mathcal{M}\tilde{\mathbf{N}}^e(\bar{\mathbf{U}}^e,\tilde{\mathbf{U}}^e))_n = -J^e\sum_{i=1}^{N_c} w^c_i G^e 
(\mathcal{I}\tilde{\nabla}l_n)_i \cdot \tilde{\mathbf{F}}((\mathcal{I}\bar{\mathbf{U}}^e)_i,(\mathcal{I}\tilde{\mathbf{U}}^e)_i) \\ + \sum_{f=0}^2 J^{ef}\sum_{j=1}^{N^f_c}w^c_j (\mathcal{I}^f l_n)_j \mathbf{n}\cdot \tilde{\mathbf{F}}^{*}((\mathcal{I}^f\bar{\mathbf{U}}^e)_j,(\mathcal{I}^f\tilde{\mathbf{U}}^e)_j),
\end{multline*}
Defining the combined differentiation and projection operator $\mathbf{P} = [\mathcal{P}_r,\mathcal{P}_s]$ via
\begin{align*}
    (\mathcal{P}_r)_{ni} &= \sum_{m=1}^{N_p}(\mathcal{M}^{-1})_{nm} w^c_i \left(\mathcal{I}\frac{\partial l_m}{\partial r}\right)_i, \\
    (\mathcal{P}_s)_{ni} &= \sum_{m=1}^{N_p}(\mathcal{M}^{-1})_{nm} w^c_i \left(\mathcal{I}\frac{\partial l_m}{\partial s}\right)_i,
\end{align*}
and the cubature lifting operators $\mathcal{L}^{f}_c$ as
\begin{equation*}
    (\mathcal{L}^{f}_c)_{nj} = \sum_{m=1}^{N_p}(\mathcal{M}^{-1})_{nm}\sum_{j=1}^{N^f_c} w^c_j (\mathcal{I}^f l_m)_j,
\end{equation*}
we can write the operator $\tilde{\mathbf{N}}^e(\bar{\mathbf{U}}^e,\tilde{\mathbf{U}}^e)$ compactly as
\begin{equation} \label{eq:KSS_3}
\tilde{\mathbf{N}}^e(\bar{\mathbf{U}}^e,\tilde{\mathbf{U}}^e) = -G^e 
\mathbf{P} \cdot \tilde{\mathbf{F}}((\mathcal{I}\bar{\mathbf{U}}^e)_i,(\mathcal{I}\tilde{\mathbf{U}}^e)_i) + \sum_{f=0}^2 \frac{J^{ef}}{J^e}\mathcal{L}^{f}_c \mathbf{n}\cdot \tilde{\mathbf{F}}^{*}((\mathcal{I}^f\bar{\mathbf{U}}^e)_j,(\mathcal{I}^f\tilde{\mathbf{U}}^e)_j).
\end{equation}
Hence, the action of the nonlinear advection can be written as the sum of the volume and surface integral contributions. The evaluation of the volume term consists of interpolating the velocity fields $\bar{\mathbf{U}}^e$ and $\tilde{\mathbf{U}}^e$ to the $N_c$ cubature nodes, followed by the actions of the combined differentiation and projection operators $\mathcal{P}_r$ and $\mathcal{P}_s$ and incorporation of the geometric factors. Similarly, the evaluation of the surface term consists of interpolating the traces of the velocity fields $\bar{\mathbf{U}}^e$ and $\tilde{\mathbf{U}}^e$ to the $N^f_c$ face cubature nodes, followed by the action of the cubature lift operator. We proceed to describe the GPU implementation and optimization of these two operations.

\subsubsection{Subcycling Advection Volume Kernel}
We show in Algorithm \ref{alg:SSV} the pseudo-code of subcycling advection volume (SAV) kernel. $N_c$ and $N_p$ threads are used for interpolation and projection steps, respectively. To perform all computations, $N_c$ threads are assigned for this kernel, unless explicitly stated otherwise.  

\begin{algorithm}[!tb]
  \caption{Subcycling Advection Volume Kernel}
  \label{alg:SSV}
\begin{boxedminipage}{1.0\textwidth}
    \begin{algorithmic}[1]
    \STATE {\bf Input:} 
     \\(1) $\mathbf{\bar{U}} = [\bar{u}, \bar{v}]$, size  $2 \times \left(K \times N_p\right)$;
    \\(2) $\mathbf{\tilde{U}} = [\tilde{u}, \tilde{v}]$, size  $2 \times \left(K \times N_p\right)$;
    \\(3) Interpolation matrix $\mathcal{I}$, size $N_c\times N_p$; 
    \\(4) Projection matrices $\mathbf{P} = [\mathcal{P}_r, \mathcal{P}_s] $, size $2\times \left( N_p\times N_c\right)$; 
    \\(5) Geometric factors $G$, size $ 4\times K$
    \STATE {\bf Output: } $\mathbf{N} = [N_u, N_v]$, size  $2 \times \left(K \times N_p\right)$;
     \FOR {$e\in\left\{1,2, \ldots K \right\}$}
     \FOR{$i\in\left\{1,2, \ldots N_c \right\}$}
    \STATE $\bar{u}_{i}=\sum_{j=1}^{N_p}\mathcal{I}_{ij}^c 
     \bar{u}_{j}^{e}$  \COMMENT{Interpolate to cubature nodes} 
    \STATE $\bar{v}_{i}=\sum_{j=1}^{N_p}\mathcal{I}_{ij}^c
    \bar{v}_{j}^{e}$ 
    \STATE $\tilde{u}_{i}=\sum_{j=1}^{N_p}\mathcal{I}_{ij}^c
     \tilde{u}_{j}^{e}$  
    \STATE $\tilde{v}_{i}=\sum_{j=1}^{N_p}\mathcal{I}_{ij}^c
     \tilde{v}_{j}^{e}$
    \STATE $F_{0;i}= \bar{u}_{i} \tilde{u}_{i}$,  $F_{1;i}= \bar{v}_{i} \tilde{u}_{i}$ \COMMENT{Compute volume flux function}    
    \STATE $F_{2;i}= \bar{u}_{i} \tilde{v}_{i}$,  $F_{3;i}= \bar{v}_{i} \tilde{v}_{i}$  
    \ENDFOR
    
    \FOR{$i\in\left\{1,2, \ldots N_c\right\}$}
    \IF[Differentiate and project back]{$i\leq N_p$}
    \STATE $r_{x} = G_{0}^{e},\; s_{x} = G_{1}^{e} $ \COMMENT{Load geometric factors}
    \STATE $r_{y} = G_{3}^{e},\; s_{y} = G_{3}^{e} $ 
    \item[] \COMMENT{Differentiate and project}
     \STATE $Fr_{0;i}= \sum_{j=1}^{N_c}\mathcal{P}_{r;ij}
    F_{0;j}$,  $Fs_{0;i}= \sum_{j=1}^{N_c}\mathcal{P}_{s;ij}
    F_{0;j}$
     \STATE $Fr_{1;i}= \sum_{j=1}^{N_c}\mathcal{P}_{r;ij}
    F_{1;j}$,  $Fs_{1;i}= \sum_{j=1}^{N_c}\mathcal{P}_{s;ij}
    F_{1;j}$
     \STATE $Fr_{2;i}= \sum_{j=1}^{N_c}\mathcal{P}_{r;ij}
    F_{2;j}$,  $Fs_{2;i}= \sum_{j=1}^{N_c}\mathcal{P}_{s;ij}
    F_{2;j}$
     \STATE $Fr_{3;i}= \sum_{j=1}^{N_c}\mathcal{P}_{r;ij}
    F_{3;j}$,  $Fs_{3;i}= \sum_{j=1}^{N_c}\mathcal{P}_{s;ij}
    F_{3;j}$
    \item[] \COMMENT{Multiply with geometric factors and update}
    \STATE $N_{u;i}^{e} = r_x Fr_{0;i} + s_x Fs_{0;i} + r_y Fr_{1;i} + s_y Fs_{1;i}$
    \STATE $N_{v;i}^{e} = r_x Fr_{2;i} + s_x Fs_{2;i} + r_y Fr_{3;i} + s_y Fs_{3;i}$
    \ENDIF
    \ENDFOR
    \ENDFOR
  \end{algorithmic}
\end{boxedminipage}
\end{algorithm}

\begin{figure}[tb]
\begin{center}
  \includegraphics[width=0.9\textwidth]{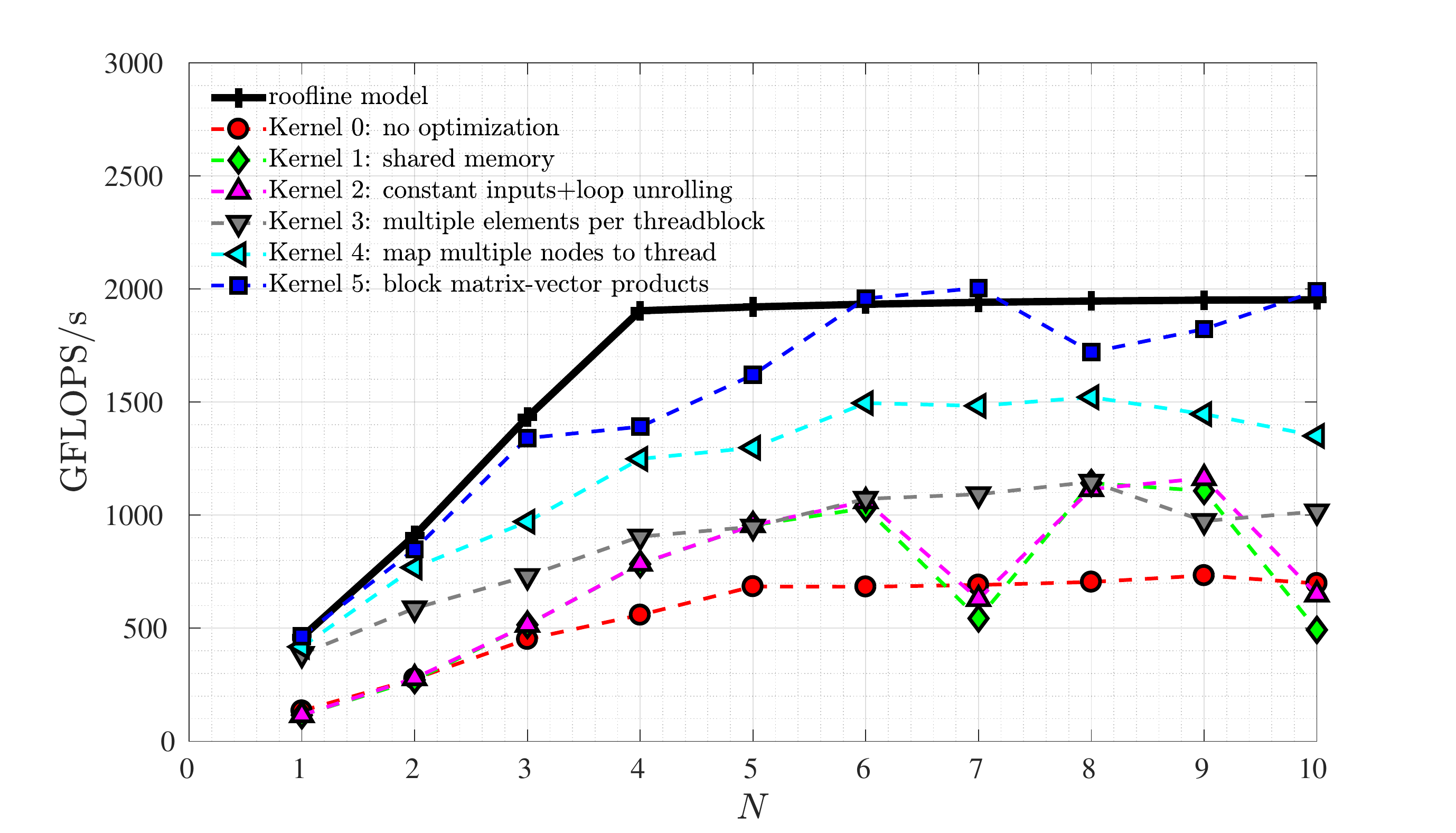}
\caption{Achieved floating point performance for the subcycling advection volume kernels compared against an empirical roofline model shown as a black line. }
\label{fig:VolumeKernel}
\end{center}
\end{figure}

We show in Figure \ref{fig:VolumeKernel} the GPU performance results of six different SAV kernels. As done above for the elliptic operator kernels, these kernels are constructed in a sequential fashion starting with a direct implementation of Algorithm \ref{alg:SSV} and applying successively optimizations. Each kernel uses the previous kernel implementation as a starting point and applies the optimizations detailed below. 

\kernelNl{SAV}{0} This kernel is a direct implementation of the pseudo-code in Algorithm \ref{alg:SSV} and serves as a reference point for measuring kernel optimizations. This kernel reads the velocity fields directly from global GPU memory during the interpolation loop stores the result in shared memory. The performance of this kernel stalls for $N\geq4$ due to excessive global memory accesses and reaches only $700$ GFLOPS/s.

\kernelNl{SAV}{1} In this kernel we introduce $4$ shared memory arrays, each with $N_p$ entries. The arrays are used to store the velocity fields before applying the interpolation operator. A memory fence is placed to ensure that all the shared memory data is loaded before the matrix-vector multiplication in the interpolation step. The reduction in global memory accesses improves the performance for $N\geq4$.  

\kernelNl{SAV}{2} In this kernel the $\texttt{const}$ qualifier is added to all unmodified input arrays, and to local variables where possible. We also label pointers with the \texttt{restrict} qualifier to explicitly state that they point to non-overlapping arrays. Additionally, all inner $\texttt{for}$ loops are unrolled, which provides the scheduler with more opportunity for instruction-level parallelism. These modifications, however, only marginally boost the performance of the kernel.     

\kernelNl{SAV}{3} In this kernel multiple elements are processed by each threadblock to better align the computational load with the hardware architecture. Running several trials, we choose the number of elements per threadblock which optimizes performance. This optimization improves the performance for low order approximations. The kernel achieves roughly 1 TFLOPS/s at high-order, which is approximately a half of the empirical shared memory bound. 

\kernelNl{SAV}{4} In this kernel, each thread processes multiple nodes of an element, in addition to each threadblock processing multiple elements. That is, each time an entry of the interpolation or projection operators is loaded from memory it can be reused multiple times in the matrix-vector multiplication. Each thread stores the interpolated variables in a register array. While this optimization yields approximately a 1.5 fold speedup, overall performance of the kernel remains lower than the shared memory bound. 

\kernelNl{SAV}{5} At high orders, the number of cubature nodes increases and becomes much larger than the number of interpolation nodes. Since each previous kernel used $N_c$ threads, as the difference between $N_c$ and $N_p$ increases most of these threads stay idle in the projection step, which reduces thread utilization and hence, negatively impacts the kernel performance. Note as well that shared memory usage for interpolated velocity fields becomes excessive with increase of interpolation orders. To avoid the thread under-utilization and the impact of shared memory latency, we use instead only $N_p$ threads with shared memory arrays of size $N_p$ for each velocity component per each element processed in the kernel. Doing so, the matrix-vector multiplication in the interpolation step is blocked and computed in multiple passes. This optimization improves the performance of the kernel substantially. The kernel achieves approximately 2 TFLOPS/s and the performance plot approaches the empirical roofline.

\subsubsection{Subcycling Advection Surface Kernel}
We show in Algorithm \ref{alg:SSS} the pseudo-code implementing the subcycling advection surface (SAS) kernel which computes the surface contribution to the subcycling advection term \eqref{eq:KSS_3}. In this kernel we require $N_f \times  N_c^f$ threads to perform the interpolation step and compute the numerical flux at the surface integration points. We then require $N_p$ thread to apply the lift operator. We therefore launch the kernel using $\max\left(N_f\times N_c^f, Np\right)$ threads per threadblock to ensure that both operations can be performed. 

\begin{algorithm}[!htb]
  \caption{Subcycling Advection Surface Kernel}
  \label{alg:SSS}
\begin{boxedminipage}{1.0\textwidth}
    \begin{algorithmic}[1]
    \STATE {\bf Input:} 
     \\(1) $\mathbf{\bar{U}} = [\bar{u}, \bar{v}]$, size  $2 \times \left(K \times N_p\right)$;
    \\(2) $\mathbf{\tilde{U}} = [\tilde{u}, \tilde{v}]$, size  $2 \times \left(K \times N_p\right)$; \\
    (3) Negative trace indices $idM$, size $K \times \left(N_f \times N_{fp}\right)$; \\
    (4) Positive trace indices $idP$, size $K \times \left(N_f \times N_{fp}\right)$; \\
    (5) Cubature Lift matrix $\mathcal{L}_c$, size $N_p \times \left(N_f \times N_c^f\right)$; \\
    (6) Interpolation matrix $\mathcal{I}^f$, size $\left(N_f\times N_c^f\right) \times N_p$; \\
    (5) Geometric factors $sG$, size $ K\times(N_f\times4)$\\
    \STATE {{\bf Output: } \\$\mathbf{N} = [N_u, N_v]$, size  $K\times N_p\times 2$}; 
     \FOR {$e\in\left\{1,2, \ldots K \right\}$}
     \FOR{$i\in\left\{1,2, \ldots \left( N_f \times N_c^f\right)\right\}$}
    \item[] \COMMENT{Interpolate to surface cubature nodes}
     \STATE $\bar{u}_{i}^{-}=\sum_{j=1}^{N_{fp}}\mathcal{I}_{ij}^c
     \bar{u}( idM^{e}_j)$,  $\bar{u}_{i}^{+}=\sum_{j=1}^{N_{fp}}\mathcal{I}_{ij}^c
     \bar{u}( idP^{e}_j)$
    \STATE $\bar{v}_{i}^{-}=\sum_{j=1}^{N_{fp}}\mathcal{I}_{ij}^c
     \bar{v}( idM^{e}_j)$,  $\bar{v}_{i}^{+}=\sum_{j=1}^{N_{fp}}\mathcal{I}_{ij}^c
     \bar{v}( idP^{e}_j)$
    \STATE $\tilde{u}_{i}^{-}=\sum_{j=1}^{N_{fp}}\mathcal{I}_{ij}^c
     \tilde{u}( idM^{e}_j)$,  $\tilde{u}_{i}^{+}=\sum_{j=1}^{N_{fp}}\mathcal{I}_{ij}^c
     \tilde{u}( idP^{e}_j)$
     \STATE $\tilde{v}_{i}^{-}=\sum_{j=1}^{N_{fp}}\mathcal{I}_{ij}^c
     \tilde{v}( idM^{e}_j)$,  $\tilde{v}_{i}^{+}=\sum_{j=1}^{N_{fp}}\mathcal{I}_{ij}^c
     \tilde{v}( idP^{e}_j)$
     \item[] \COMMENT{Compute flux function}
    \STATE $n_{x} = sG_{0}^{ef},\; n_{y} = sG_{1}^{ef}  \;J^{ef} = sG_{2}^{ef},\; (J^{e})^{-1} = sG_{3}^{e} $ 
    \STATE $\alpha = 0.5 (J^{e})^{-1} J^{ef}$, $\lambda_i = \max(|n_x\bar{u}_{i}^{-}+n_y\bar{v}_{i}^{-}|,|n_x\bar{u}_{i}^{+}+n_y\bar{v}_{i}^{+}|)$
    \STATE $\mathbf{F}^*_{u;i}=\alpha\left(n_x \left(\bar{u}_{i}^{+}\tilde{u}_{i}^{+}+\bar{u}_{i}^{-}\tilde{u}_{i}^{-} \right) 
    +n_y\left( \bar{v}_{i}^{+}\tilde{u}_{i}^{+}+\bar{v}_{i}^{-}\tilde{u}_{i}^{-}\right) + \lambda_{i}\left( \tilde{u}_{i}^{-}-\tilde{u}_{i}^{+} \right) 
     \right)$
     \STATE $\mathbf{F}^*_{v;i}= \alpha\left( n_x \left(\bar{u}_{i}^{+}\tilde{v}_{i}^{+}+\bar{u}_{i}^{-}\tilde{v}_{i}^{-} \right) 
    +n_y\left( \bar{v}_{i}^{+}\tilde{v}_{i}^{+}+\bar{v}_{i}^{-}\tilde{v}_{i}^{-}\right) +
    \lambda_{i}\left( \tilde{v}_{i}^{-}-\tilde{v}_{i}^{+} \right) 
     \right)$
    \ENDFOR
    \FOR{$i\in\left\{1,2, \ldots N_p\right\}$}
    \STATE $N_{u;i}=\sum_{j=1}^{N_f \times N_{c}^f}\mathcal{L}_{c;ij} \mathbf{F}^*_{u;i}$    \COMMENT{Lift numerical flux}   
    \STATE $N_{v;i}=\sum_{j=1}^{N_f \times N_{c}^f}\mathcal{L}_{c;ij} \mathbf{F}^*_{v;i}$ 
    \STATE $N_{u;i}^{(e)}  += N_{u;i}$ \COMMENT{Add to volume contribution}
    \STATE $N_{v;i}^{(e)} += N_{v;i}$
    \ENDFOR
    \ENDFOR
  \end{algorithmic}
\end{boxedminipage}
\end{algorithm}

\begin{figure}[htb]
\begin{center}
  \includegraphics[width=0.9\textwidth]{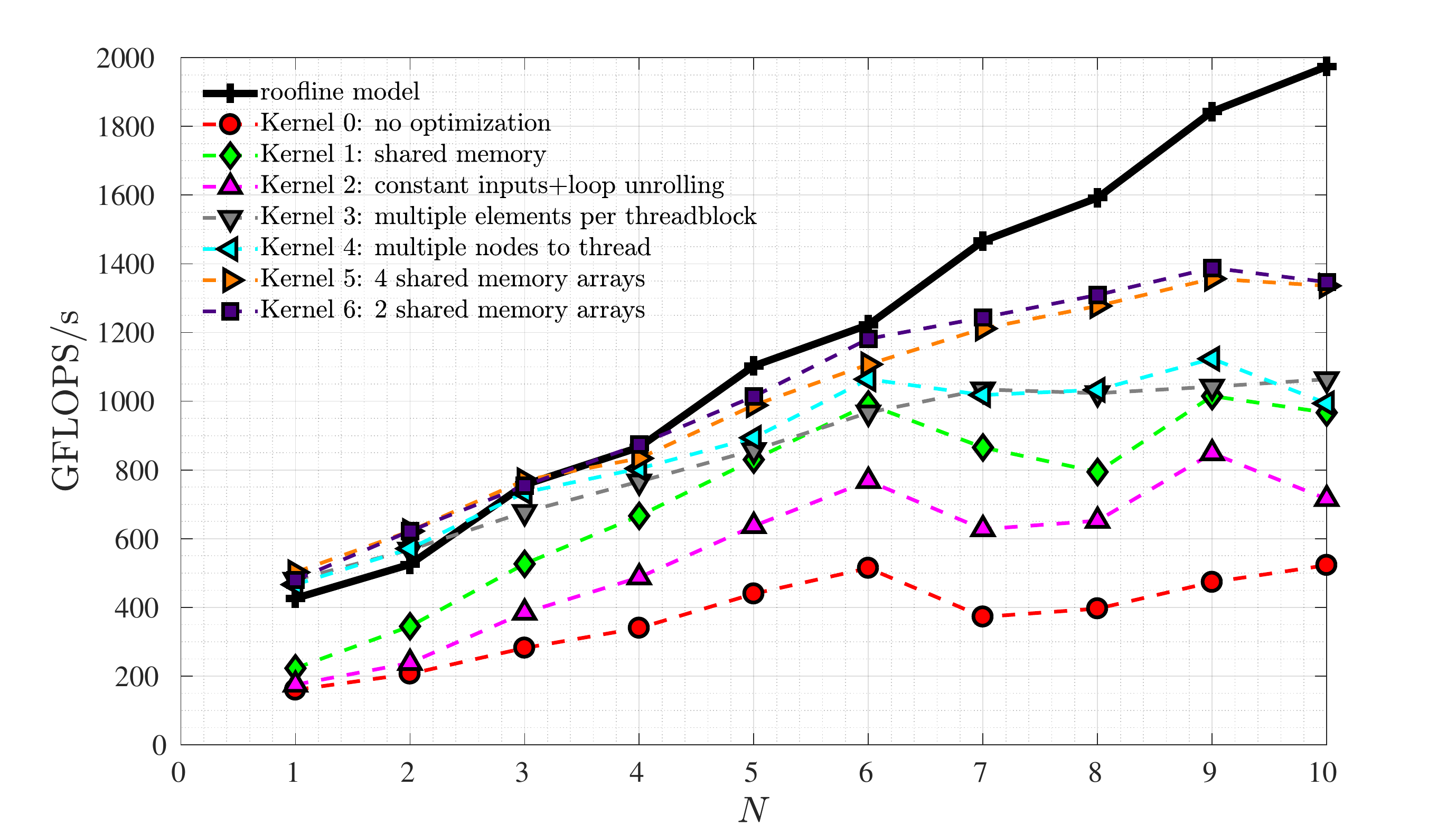}
\caption{Achieved floating point performance for the subcycling surface kernels compared against an empirical roofline model shown as a black line.}
\label{fig:SurfaceKernel}
\end{center}
\end{figure}

We show in Figure \ref{fig:SurfaceKernel} the GPU performance of seven separate kernels implemented to compute the surface contribution to the subcycling advection term. As described above for previous kernels, these kernels are constructed using sequential optimization steps, starting from the direct implementation of Algorithm \ref{alg:SSS}. We detail the optmizations performed in each kernel below.

\kernelNl{SAS}{0} This kernel is a direct implementation of the pseudo-code in Algorithm \ref{alg:SSS} and serves as a reference point for measuring kernel optimizations. This kernel uses two shared memory arrays of size $N_f\times N_c^f$ to store the numerical flux for surface integration points. Each of the velocity fields are loaded directly from the global memory in the interpolation step. The excessive global memory accesses limit the performance of this kernel and performance reaches only $400$ GFLOPS/s, which is one fifth of the predicted empirical roofline for $N=10$. 

\kernelNl{SAS}{1} In this kernel we introduce eight additional shared memory arrays of size $N_f \times N_{fp}$ to store the internal and neighbour trace data of the velocity fields. All the required data is loaded from global memory at the beginning of the kernel, before the interpolation step. The resulting reduction in global memory reads significantly improves the performance of the kernel and performance $800$ GFLOPS/s, which is a two-fold speedup compared with SAS Kernel 0.  

\kernelNl{SAS}{2} In this kernel we add the \texttt{const} qualifier to all unmodified input variables. We also label input pointers with the \texttt{restrict} qualifier to explicitly state that they point to non-overlapping arrays. Additionally, all serial $\texttt{for}$ loops in interpolation and lifting steps are unrolled to increase instruction-level parallelism. Although these modifications provide further optimization opportunities for the compiler, our results indicate that they have only a minor effect on the achieved performance. 

\kernelNl{SAS}{3} In this kernel multiple elements are processed by each threadblock to better balance the occupancy and the data movement. As for the volume kernel, the optimal number of elements per threadblock is optimized by testing over several options. The performance improvement resulting from this optimization is modest, and much better at low-order approximations. This kernel performs around $1$ TFLOPS/s for $N>5$ which is $50 \%$ of the empirical bound for $N=10$.   

\kernelNl{SAS}{4} In this kernel multiple nodes of different elements are processed by a thread to further increase the occupancy and to reuse fetched interpolation and lift operators. This optimization slightly improves performance at low order approximations. However, due to excessive shared memory requirements we cannot load a sufficient number of elements in a single thread block to make this optmization yield a performance improvement at high orders.

\kernelNl{SAS}{5} In this kernel shared memory usage is reduced by a factor of two. We first load the  velocity fields from global memory to shared memory arrays and then interpolate the surface integration points. The interpolated velocity fields are stored in register arrays and loaded back to the same shared memory arrays after local memory barrier. This reduction in shared memory usage allows us to load more elements per thread block and take advantage of the optmizations performed in the previous kernel giving  an approximate $20\%$ performance improvement for $N>5$. The kernel reaches $1.35$ TFLOPS/s.   

\kernelNl{SAS}{6} In this kernel shared memory usage is further reduced by a factor of two using two additional thread synchronizations. This kernel utilizes only two shared memory arrays where velocity components are loaded and interpolated to the integration nodes in sets of two before each thread synchronization. We process one velocity field by fetching interior and exterior trace values from the global memory to increase the likelihood of data caching. Performance is slightly improved achieving $1.4$ TFLOPS/s.

\section{Conclusion}
In this study, we presented a GPU-optimized high-order discontinuous Galerkin method for approximating the incompressible Navier-Stokes equations. To reduce the cost of each semi-implicit time step we use a semi-Lagrangian subcycling approach. Performance studies show that this approach shifts the computational load away from the linear solvers towards the explicit advection stage. 

We presented an empirical performance roofline model to assist in quantifying GPU performance as well as indicate when kernels are performing near empirical limits. We conducted a detailed study of the most computationally intensive kernels in the linear solver stage as well as the subcycling advection stage. We detailed the optimization of each of these kernels targeting the Nvidia Tesla P100 GPU. The resulting performance measures of the optimized kernels indicate that the solver is performing well on the GPU.

The GPU performance of the three dimensional versions of each of the high-order operators in the INS scheme for tetrahedral elements remains to be investigated. Furthermore, significant performance gains can potentially be obtained by considering modifications to aspects of the scheme such as more sophisticated preconditioning techniques and polynomial bases which sparsity finite-element operators. These topics will studied in future works.


\section{Acknowledgements}
This research was supported in part by the Exascale Computing Project (17-SC-20-SC), a collaborative effort of two U.S. Department of Energy organizations (Office of Science and the National Nuclear Security Administration) responsible for the planning and preparation of a capable exascale ecosystem, including software, applications, hardware, advanced system engineering, and early testbed platforms, in support of the nations exascale computing imperative.

In addition, the authors would like to kindly acknowledge Advance Research Computing at Virginia Tech for providing readily accessible computational resources. Finally, this research was supported in part by the John K. Costain Faculty Chair in Science at Virginia Tech.

\printbibliography
\end{document}